\documentclass[a4paper, final, 12pt]{article}
\usepackage{amsmath, amssymb, latexsym, amscd, amsthm,amsfonts,amstext}
\usepackage[mathscr]{eucal}
\usepackage{graphicx}
\usepackage{subfig}
\usepackage{float}
\usepackage{color}
\usepackage{hyperref}
\usepackage[utf8]{inputenc}
\usepackage[english]{babel}

\numberwithin{equation}{section}
\input{tcilatex}
\begin{document}

\title{ On the Mean Field Games System With Lateral Cauchy Data via Carleman
Estimates\thanks{\textbf{Funding}. The work of J. Li was partially supported
by the NSF of China No. 11971221, Guangdong NSF Major Fund No. 2021ZDZX1001,
the Shenzhen Sci-Tech Fund No. RCJC20200714114556020, JCYJ20200109115422828
and JCYJ20190809150413261, National Center for Applied Mathematics Shenzhen,
and SUSTech International Center for Mathematics. The work of H. Liu was
supported by the Hong Kong RGC General Research Funds (projects 12302919,
12301420 and 11300821) and the France-Hong Kong ANR/RGC Joint Research
Grant, A-CityU203/19.}}
\author{ Michael V. Klibanov \thanks{
Department of Mathematics and Statistics, University of North Carolina at
Charlotte, Charlotte, NC, 28223, USA, mklibanv@uncc.edu} \and Jingzhi Li 
\thanks{
Department of Mathematics \& National Center for Applied Mathematics
Shenzhen \& SUSTech International Center for Mathematics, Southern
University of Science and Technology, Shenzhen 518055, P.~R.~China,
li.jz@sustech.edu.cn} \and Hongyu Liu \thanks{
Department of Mathematics, City University of Hong Kong, Kowloon, Hong Kong
SAR, P.R. China, hongyliu@cityu.edu.hk} }
\date{}
\maketitle

\begin{abstract}
The second order Mean Field Games system (MFGS) in a bounded domain with the
lateral Cauchy data is considered. This means that both Dirichlet and
Neumann boundary data for the solution the MFGS are given. Two H\"{o}lder
stability estimates for two slightly diffeent cases are derived. These
estimates indicate how stable the solution of the MFGS is with respect to
the possible noise in the lateral Cauchy data. Our stability estimates imply
uniqueness. The key mathematical apparatus is the apparatus of two new
Carleman estimates.
\end{abstract}

\textbf{Key Words}: mean field games system, lateral Cauchy data, ill-posed
and inverse problems, Carleman estimates, H\"{o}lder stability, uniqueness

\textbf{2020 MSC codes}: 35R30, 91A16

\section{Introduction}

\label{sec:1}

The Mean Field Games (MFG) theory studies the behavior of infinitely many
agents/players, whose intent is to optimize their values depending on the
dynamics of the whole set of players. This theory is enjoying a rapidly
growing number of applications to various socio-economic processes. Some
applications are in, e.g. economics \cite{A,A2,Burger,LL2,Trusov}, election
dynamics \cite{Chow}, combating corruption \cite{KB}-\cite{Kol}, crowd
motion \cite{A1}, sociology \cite{B} and interaction of electric vehicles 
\cite{Cou}.

It is pointed out in \cite{Burger} that the most elegant feature of the MFG
theory is that it is governed by a system of PDEs, the so-called Mean Field
Games System (MFGS), see, e.g. \cite{A,LL1,LL2,LL3} for this system. Even
though different specific applications are governed by different versions of
the MFGS, all these versions have similar features. Hence, given a rapidly
growing role of socio-economic sciences in the modern society, it is
important to study various mathematical questions for at least one MFGS, and
this is what we do in the current paper.

The MFG theory was first introduced in the seminal publications of Lasry and
Lions \cite{LL1,LL2,LL3} and Huang, Caines and Malham\'{e} \cite%
{Huang1,Huang2}. In these works, a second order system of two coupled
nonlinear parabolic equations was derived. Let $x\in \mathbb{R}^{n}$ be the
state variable characterizing position of an agent and $t>0$ be time. These
two equations have two opposite directions of time and they form the second
order Mean Field Games System. The first equation is the
Hamilton-Jacobi-Bellman equation (HJB) for the value function $u\left(
x,t\right) $ and the second equation is the Fokker-Planck equation (FP) for
the function $m\left( x,t\right) ,$ which describes the density of agents.

We assume that $x\in \Omega ,$ where $\Omega \subset \mathbb{R}^{n}$ is a
bounded domain with the piecewise smooth boundary $\partial \Omega .$ We
call a problem for the MFGS \textquotedblleft forward problem" if it is
required to find a solution of the MFGS, assuming that all coefficients of
this system are known. And we call a problem for the MFGS \textquotedblleft
coefficient inverse problem" if it is required to reconstruct both an
unknown coefficient(s) of this system and its solution, using boundary
measurements.

It is assumed in the traditional setting of forward problems for the MFGS
that functions 
\begin{equation}
u\left( x,T\right) ,\text{ }m\left( x,0\right) ,\text{ }x\in \Omega
\label{1.1}
\end{equation}%
are known, which are terminal and initial conditions for functions $u\left(
x,t\right) $ and $m\left( x,t\right) .$ Uniqueness in this case is
established only under a very restrictive the so-called \textquotedblleft
monotonicity" condition \cite{LL3}. On the other hand, Klibanov and
Averboukh were the first ones who have introduced Carleman estimates in the
MFG theory \cite{MFG1}. This allowed them to obtain Lipschitz stability
estimate for the solution of the forward problem for the MFGS in the case
when, instead of (\ref{1.1}), the data are: 
\begin{equation}
u\left( x,T\right) ,\text{ }m\left( x,0\right) ,\text{ }m\left( x,T\right) ,%
\text{ }x\in \Omega .  \label{1.2}
\end{equation}%
In \cite{MFG2} Carleman estimates were applied to obtain Lipschitz stability
estimate for the case when (\ref{1.2}) is replaced with 
\begin{equation}
u\left( x,T\right) ,\text{ }u\left( x,0\right) ,\text{ }m\left( x,T\right) ,%
\text{ }x\in \Omega .  \label{1.3}
\end{equation}%
Finally, in \cite{MFG4} H\"{o}lder stability estimates for two forward
problems for the MFGS are obtained in the case when three functions in
either (\ref{1.2}) or (\ref{1.3}) are replaced with either of two pairs:%
\begin{equation}
\left( u\left( x,T\right) ,\text{ }m\left( x,T\right) \right) \text{ or }%
\left( u\left( x,0\right) ,\text{ }m\left( x,0\right) \right) .  \label{1.4}
\end{equation}%
Uniqueness of corresponding forward problems follows immediately from those
stability estimates. An additional importance of stability estimates is due
to the fact that they provide \emph{a priori }estimates of the stability of
solutions with respect to a possible noise in the input data. Note that the
input data are always noisy in applications.

In all cases (\ref{1.1})-(\ref{1.4}) the data are either at $\left\{
t=0\right\} $ or at $\left\{ t=T\right\} $ or at both. The goal of this
paper is to address the following question: \emph{Would it be possible to
obtain some stability and uniqueness results for the MFGS if these data
would be replaced with Dirichlet and Neumann boundary conditions at the
lateral boundary, which is }$\partial \Omega \times \left( 0,T\right) $\emph{%
? }

Our answer to this question is positive: we obtain H\"{o}lder stability
estimates for the MFGS in the case when the lateral Cauchy data are given,
i.e. both Dirichlet and Neumann boundary conditions are known. These
estimates immediately imply uniqueness of the solution of that problem. The
main tool here is again the tool of Carleman estimates.

Recently H\"{o}lder and Lipschitz stability estimates were obtained in \cite%
{MFG5,MFG6} for some coefficient inverse problems for the MFGS using the
framework of the paper \cite{BukhKlib}. In \cite{BukhKlib} Carleman
estimates were introduced in the field of coefficient inverse problems for
the first time, also, see, e.g. \cite{Ksurvey,Yam} and the book \cite{KL}
for this method. In \cite{MFG7} the forward problem of \cite{MFG1} with the
data (\ref{1.2}) is solved numerically via a globally convergent so-called
convexification method. In \cite{MFG8} a coefficient inverse problem for the
MFGS is numerically solved by another version of the convexification method.
Convexification is also based on Carleman estimates. See, e.g. the book \cite%
{KL} for a detailed description of both the technique of \cite{BukhKlib} and
the convexification method.

While the above cited works \cite{MFG1}-\cite{MFG8} study the MFGS for the
case of the data generated by a single measurement event, publications \cite%
{Liu1,Liu2,Liu3} are concerned with proofs of uniqueness theorems for the
MFGS for the case of the data generated by multiple measurements.

There are three main difficulties of the MFGS, as compared with a single
parabolic PDE:

\begin{enumerate}
\item Two PDEs forming the MFGS have opposite directions of time. Hence, the
classical theory of parabolic PDEs is inapplicable to the MFGS.

\item The presence of the integral operator 
\begin{equation}
\dint\limits_{\Omega }K\left( x,y\right) m\left( y,t\right) dy  \label{1.5}
\end{equation}%
in the HJB equation. This is a very unusual term characterizing the global
interaction of players.

\item The presence of the Laplace operator with $\Delta u\left( x,t\right) $
in the FP equation for the function $m\left( x,t\right) .$Indeed, the term $%
\Delta u\left( x,t\right) $ is involved in the principal part of the HJB
equation.
\end{enumerate}

Summarizing, the combination of items 1-3 does not allow one to
automatically \textquotedblleft project" the methodology developed for
ill-posed problems and CIPs for a single parabolic equation on corresponding
problems for the MFGS, see \cite[section 2.3]{MFG6} for a similar conclusion.

\textbf{Remarks 1.1}:

\begin{enumerate}
\item \emph{The minimal smoothness assumptions are of a low priority in the
theory of Ill-Posed and Inverse Problems, see, e.g. \cite{Nov1,Nov2}, \cite[%
Theorem 4.1]{Rom}. Similarly, we are not concerned about these assumptions
in the current paper. Indeed, nonlinear ill-posed and inverse problems are
quite challenging ones in their own rights. Thus, it makes sense to obtain,
at least on the first stage of the research, results for sufficiently smooth
functions. And then to decrease the required smoothness on the second stage
of research. }

\item \emph{A similar statement is true with the reference to our assumption
in section 2 that the domain of our interest }$\Omega $\emph{\ is a
rectangular prism. We use this assumption on the first stage of our research
of the problem of this paper: to simplify the presentation. More general
domains might be considered in the future research in this direction. }
\end{enumerate}

This paper is organized as follows. In section 2 we state the problem of our
interest here. In section 3 we present two forms of the kernel $K\left(
x,y\right) $ in (\ref{1.5}) we work with. In section 4 we formulate two
Carleman estimates we work with and prove one of them. The proof of the
second one is very similar and is, therefore, omitted here. In section 5 we
prove our main results proclaimed above: two H\"{o}lder stability estimates
for two slightly different cases. All functions below are real valued ones.

\section{Problem Statement}

\label{sec:2}

When initial and terminal data (\ref{1.2})-(\ref{1.4}) are given, as in \cite%
{MFG1}-\cite{MFG4}, we can work with an almost arbitrary kernel $K\left(
x,y\right) $ of the integral operator (\ref{1.5}). However, since we work
here with lateral Cauchy data, then possible forms of kernels are limited,
see section 3. These limitations, in turn lead to a limitation of a possible
form of the Carleman Weight Function, which should depend then only on $%
x_{1} $ and $t$ (subsection 3.3). This, in turn leads to a limitation of the
shape of the domain $\Omega $ we work with.

To single out the variable $x_{1},$we set $x=\left( x_{1},\overline{x}%
\right) ,$ $\overline{x}=\left( x_{2},...,x_{n}\right) .$ Let $a,b,a_{i},T>0$
be some real numbers, where $a<b$ and $i=2,...,n.$ Define the rectangular
prism $\Omega \subset \mathbb{R}^{n}$ with its boundary $\partial \Omega $
as well as associated sets as:%
\begin{equation}
\Omega =\left\{ x:a<x_{1}<b,-a_{i}<x_{i}<a_{i},i=2,...,n\right\} ,
\label{2.1}
\end{equation}%
\begin{equation}
\Omega _{1}=\left\{ \overline{x}:-a_{i}<x_{i}<a_{i},i=2,...,n\right\} ,
\label{2.01}
\end{equation}%
\begin{equation}
\partial _{+}^{1}\Omega =\left\{ x\in \partial \Omega :x_{1}=b\right\}
,\partial _{+}^{1}\Omega _{T}=\partial _{+}^{1}\Omega \times \left(
0,T\right) ,  \label{2.2}
\end{equation}%
\begin{equation}
\partial _{-}^{1}\Omega =\left\{ x\in \partial \Omega :x_{1}=a\right\}
,\partial _{-}^{1}\Omega _{T}=\partial _{-}^{1}\Omega \times \left(
0,T\right) ,  \label{2.3}
\end{equation}%
\begin{equation}
\partial _{\pm }^{i}\Omega =\left\{ x\in \partial \Omega :x_{i}=\pm
a_{i}\right\} ,\text{ }\partial _{\pm }^{i}\Omega _{T}=\partial _{\pm
}^{i}\Omega \times \left( 0,T\right) ,\text{ }i=2,..,n,  \label{2.4}
\end{equation}%
\begin{equation}
Q_{T}=\Omega \times \left( 0,T\right) ,S_{T}=\partial \Omega \times \left(
0,T\right) =\cup _{i=1}^{n}\left( \partial _{\pm }^{i}\Omega _{T}\right) .
\label{2.5}
\end{equation}%
Since $\partial \Omega $ is not smooth, then, using (\ref{2.1})-(\ref{2.5}),
we need to define spaces $H^{2,1}\left( S_{T}\right) $ and $H^{1,0}\left(
S_{T}\right) ,$%
\begin{equation}
\left. 
\begin{array}{c}
H^{2,1}\left( S_{T}\right) =\left\{ p:\left\Vert p\right\Vert
_{H^{2,1}\left( S_{T}\right) }^{2}=\dsum\limits_{i=1}^{n}\left\Vert
p\right\Vert _{H^{2,1}\left( \partial _{\pm }^{i}\Omega _{T}\right)
}^{2}<\infty \right\} , \\ 
H^{1,0}\left( S_{T}\right) =\left\{ p:\left\Vert p\right\Vert
_{H^{1,0}\left( S_{T}\right) }^{2}=\dsum\limits_{i=1}^{n}\left\Vert
p\right\Vert _{H^{1,0}\left( \partial _{\pm }^{i}\Omega _{T}\right)
}^{2}<\infty \right\} .%
\end{array}%
\right.  \label{2.50}
\end{equation}

Let $\alpha =const.>0.$ The MFGS of our study is:%
\begin{equation}
\left. 
\begin{array}{c}
u_{t}(x,t)+\alpha \Delta u(x,t)-k(x)(\nabla u(x,t))^{2}/2+ \\ 
+P\left( x,t,\dint\limits_{\Omega }K\left( x,y\right) m\left( y,t\right)
dy,m\left( x,t\right) \right) =0,\text{ }\left( x,t\right) \in Q_{T}, \\ 
m_{t}(x,t)-\alpha \Delta m(x,t)-\func{div}(k(x)m(x,t)\nabla u(x,t))=0,\text{ 
}\left( x,t\right) \in Q_{T},%
\end{array}%
\right.  \label{2.6}
\end{equation}%
see, e.g. \cite{A}. The function $k(x)$ in (\ref{2.6}) characterizes the
reaction of the controlled agent to an action applied at the point $x$ \cite[%
section 5]{MFG1}. The function $P$ in (\ref{2.6}) is the interaction term.
The integral operator in it is called the \textquotedblleft global
interaction term". This term expresses the average action on an agent, who
is located at the point $x,$ by all other agents \cite[section 5]{MFG1}. The
case with $K\left( x,y\right) \equiv 0$ is called \textquotedblleft local
interaction" and has a lesser applied meaning.

As stated in section 1, we assume that both Dirichlet and Neumann boundary
conditions at $S_{T}$ are known, i.e. the lateral Cauchy data are given. Let 
$\nu =\nu \left( x\right) $ be the unit outward looking normal vector at
such a point $\left( x,t\right) \in S_{T},$ which belongs to one of smooth
pieces $\partial _{\pm }^{i}\Omega _{T},$ see (\ref{2.4}), (\ref{2.5}).
Thus, we assume that the following functions $g_{0},g_{1},p_{0},p_{1}$ are
known:%
\begin{equation}
u\mid _{S_{T}}=g_{0}\left( x,t\right) ,\text{ }\partial _{\nu }u\mid
_{S_{T}}=g_{1}\left( x,t\right) ,  \label{2.8}
\end{equation}%
\begin{equation}
m\mid _{S_{T}}=p_{0}\left( x,t\right) ,\text{ }\partial _{\nu }p\mid
_{S_{T}}=p_{1}\left( x,t\right) .  \label{2.9}
\end{equation}

Denote 
\begin{equation}
Q_{T,k\varepsilon }=\Omega \times \left( k\varepsilon ,T-k\varepsilon
\right) \subset Q_{T},\text{ }\forall \varepsilon \in \left( 0,\frac{T}{4}%
\right) ,\text{ }k=1,2.  \label{3.01}
\end{equation}

\textbf{Problem Statement}. \emph{Assume that there exist two 6-d vector
functions satisfying (\ref{2.6})-(\ref{2.9}):}%
\begin{equation}
\left. 
\begin{array}{c}
\left( u_{1},m_{1},g_{0,1},g_{1,1},p_{0,1},p_{1,1}\right) \text{ and }\left(
u_{2},m_{2},g_{0,2},g_{1,2},p_{0,2},p_{1,2}\right) , \\ 
u_{1},m_{1},u_{2},m_{2}\in H^{3}\left( Q_{T}\right) .%
\end{array}%
\right.  \label{2.13}
\end{equation}%
\emph{For an arbitrary }$\varepsilon \in \left( 0,T/3\right) $ \emph{%
estimate norms }$\left\Vert u_{1}-u_{2}\right\Vert _{H^{2,1}\left(
Q_{T,\varepsilon }\right) },$\emph{\ }$\left\Vert m_{1}-m_{2}\right\Vert
_{H^{2,1}\left( Q_{T,\varepsilon }\right) }$\emph{\ via norms of the
differences of the lateral Cauchy data: }%
\begin{equation}
\left\Vert g_{0,1}-g_{0,2}\right\Vert _{H^{2,1}\left( S_{T}\right)
},\left\Vert g_{1,1}-g_{1,2}\right\Vert _{H^{1,0}\left( S_{T}\right) },
\label{2.15}
\end{equation}%
\begin{equation}
\left\Vert p_{0,1}-p_{0,2}\right\Vert _{H^{2,1}\left( S_{T}\right)
},\left\Vert p_{1,1}-p_{1,2}\right\Vert _{H^{1,0}\left( S_{T}\right) }.
\label{2.16}
\end{equation}

\section{Two Forms of the Kernel $K\left( x,y\right) $}

\label{sec:3}

Just as in \cite{MFG6}, we work in this paper with two forms of the kernel $%
K\left( x,y\right) $ of the integral operator in the first equation (\ref%
{2.6}).

\subsection{The first form}

\label{sec:3.1}

Let $\delta \left( z\right) ,z\in \mathbb{R}$ be the delta function and $%
R_{1}>0$ be a number. We assume that the kernel of the function $K\left(
x,y\right) $ in (\ref{1.5}) and (\ref{2.6}) has the following form:%
\begin{equation}
K\left( x,y\right) =\delta \left( y_{1}-x_{1}\right) K_{1}\left( x,\overline{%
y}\right) ,  \label{30.02}
\end{equation}%
\begin{equation}
K_{1}\in L_{\infty }\left( \Omega \times \Omega _{1}\right) ,\left\Vert
K_{1}\right\Vert _{L_{\infty }\left( \Omega \times \Omega _{1}\right) }\leq
R_{1},  \label{30.03}
\end{equation}%
where the domain $\Omega _{1}$ is defined in (\ref{2.01}). By (\ref{30.02}) 
\begin{equation}
\dint\limits_{\Omega }K\left( x,y\right) m\left( y,t\right)
dy=\dint\limits_{\Omega _{1}}K_{1}\left( x,\overline{y}\right) m\left( x_{1},%
\overline{y},t\right) d\overline{y}.  \label{30.04}
\end{equation}

\subsection{The second form}

\label{sec:3.2}

For $z\in \mathbb{R},$ let $H\left( z\right) $ be the Heaviside function,%
\begin{equation*}
H\left( z\right) =\left\{ 
\begin{array}{c}
1,\text{ }z>0, \\ 
0,\text{ }z<0.%
\end{array}%
\right.
\end{equation*}%
The second form of the function $K\left( x,y\right) $ we work with is:%
\begin{equation}
K\left( x,y\right) =H\left( y_{1}-x_{1}\right) K_{2}\left( x,y\right) ,
\label{30.05}
\end{equation}%
\begin{equation}
K_{2}\in L_{\infty }\left( \Omega \times \Omega \right) ,\left\Vert
K_{2}\right\Vert _{L_{\infty }\left( \Omega \times \Omega \right) }\leq
R_{1}.  \label{30.06}
\end{equation}%
It follows from (\ref{2.1}) and (\ref{2.01}) that the integral in (\ref{2.6}%
) becomes in this case:%
\begin{equation}
\dint\limits_{\Omega }K\left( x,y\right) m\left( y,t\right)
dy=\dint\limits_{x_{1}}^{b}\left( \dint\limits_{\Omega _{1}}K_{2}\left(
x,y_{1},\overline{y}\right) m\left( y,t\right) d\overline{y}\right) dy_{1}.
\label{30.07}
\end{equation}

\subsection{The Carleman Weight Function and estimates of the integral
operators in (\protect\ref{30.04}) and (\protect\ref{30.07})}

\label{sec:3.3}

Let $c>0$ be a number and $\lambda >1$ be a large parameter. Introduce the
Carleman Weight Function (CWF) $\psi _{\lambda }\left( x_{1},t\right) ,$%
\begin{equation}
\psi _{\lambda }\left( x_{1},t\right) =\exp \left[ \lambda \left(
x_{1}^{2}-c^{2}\left( t-T/2\right) ^{2}\right) \right] .  \label{3.1}
\end{equation}%
By (\ref{2.1})%
\begin{equation}
\exp \left[ 2\lambda \left( a^{2}-c^{2}\left( t-T/2\right) ^{2}\right) %
\right] \leq \psi _{\lambda }^{2}\left( x_{1},t\right) \leq e^{2\lambda
b^{2}}\text{ in }Q_{T}.  \label{3.100}
\end{equation}%
The unusual feature of the CWF (\ref{3.1}) is that its dependence from only
one spatial variable is combined with the presence of only a single large
parameter $\lambda $ rather than two large parameters, as it is in the
conventional cases of CWFs for parabolic operators, see, e.g. \cite{Ksurvey}%
, \cite[section 2.3]{KL}, \cite[\S 1 of chapter 4]{LRS}. The reason of the
dependence of $\psi _{\lambda }\left( x_{1},t\right) $ from only one spatial
variable is that we consider only two above types of the kernel $K\left(
x,y\right) $ in the integral operator in (\ref{1.5}) and (\ref{2.6}). The
reason of the absence of the second large parameter in $\psi _{\lambda
}\left( x_{1},t\right) $ is our desire to study in the future the problem of
this paper numerically, using the convexification numerical method, see,
e.g. \cite{KLpar,MFG8} and \cite[chapter 9]{KL} for applications of this
method to a CIP for a single parabolic PDE, and also see \cite{MFG8} for the
case of the MFGS \cite{MFG8}. The key role in the convexification is played
by a CWF. However, if a CWF depends on two large parameters rather than just
one, then it changes too rapidly, which is inconvenient for the numerical
implementation.

Lemma 3.1 follows immediately from Lemmata 3.1 and 3.2 of \cite{MFG6}.

\textbf{Lemma 3.1.} \emph{Let }$\psi _{\lambda }\left( x_{1},t\right) $\emph{%
\ be the CWF defined in (\ref{3.1}). Let }$K\left( x,y\right) $\emph{\ be
any of two functions (\ref{30.02}) or (\ref{30.05}). In the case of (\ref%
{30.02}) the integral in (\ref{1.5}) and (\ref{2.6}) is understood as in (%
\ref{30.04}). Let (\ref{30.03}) and (\ref{30.06}) hold. Then the following
estimate is valid:}%
\begin{equation*}
\dint\limits_{Q_{T}}\left( \dint\limits_{\Omega }K\left( x,y\right) f\left(
y,t\right) dy\right) ^{2}\psi _{\lambda }^{2}dxdt\leq
B\dint\limits_{Q_{T}}f^{2}\left( x,t\right) \psi _{\lambda }^{2}\left(
x,t\right) dxdt,\text{ }\forall f\in L_{2}\left( Q_{T}\right) ,\text{ }%
\forall \lambda >0,
\end{equation*}%
\emph{where the number }$B=B\left( R,Q_{T}\right) >0$\emph{\ depends only on
listed parameters.}

\section{Carleman Estimates}

\label{sec:4}

CWF (\ref{3.1}) was proved in the past only in \cite{KLpar} and \cite[%
chapter 9]{KL}, where an analog of the Carleman estimate of Theorem 4.1 was
obtained. However, it was assumed in \cite{KLpar,KL} that the test function $%
u\left( x,t\right) $ has zero Dirichlet boundary condition at the entire
surface $S_{T}$ as well as zero Neumann boundary condition at $\partial
_{+}^{1}\Omega _{T}\subset S_{T},$ see (\ref{2.1})-(\ref{2.2}). On the other
hand, we cannot assume that these boundary conditions equal zero since we
are supposed to provide estimates via them, see (\ref{2.15}), (\ref{2.16}).
Therefore, it is necessary to prove Theorem 4.1.

To obtain our target H\"{o}lder stability estimate via boundary functions (%
\ref{2.8}), (\ref{2.9}), we need to carefully work with terms under $\left(
\cdot \right) _{x_{i}}$ and $\left( \cdot \right) _{t}$ signs in our
pointwise Carleman estimates, and this is what we do below in the proof of
Theorem 4.1. We also remind that derivations of Carleman estimates are
always space consuming, see, e.g. \cite{KLpar,KL,LRS,Yam}. A Carleman
estimate needs to be derived only for the principal part of a PDE operator
since it is independent on low order derivatives of this operator \cite[%
Lemma 2.1.1]{KL}. In Theorem 4.1, we slightly abuse the notation for $%
u\left( x,t\right) $ in (\ref{2.6}).

To save space, we formulate two Carleman estimates at once: for operators $%
\partial _{t}-\alpha \Delta $ and $\partial _{t}+\alpha \Delta .$ Their
proofs are almost identical for both cases.

\textbf{Theorem 4.1.}\emph{\ Let conditions (\ref{2.1})-(\ref{2.5}) and (\ref%
{3.1}) hold. Then there exists a sufficiently large number }$\lambda
_{0}=\lambda _{0}\left( \alpha ,c,Q_{T}\right) \geq 1$\emph{\ and a number }$%
C=C\left( \alpha ,c,Q_{T}\right) >0,$\emph{\ both depending only on listed
parameters, such that the following Carleman estimate holds for all }$%
\lambda \geq \lambda _{0}$\emph{\ and all functions }$u\in H^{3}\left(
Q_{T}\right) :$ 
\begin{equation*}
\dint\limits_{Q_{T}}\left( u_{t}\pm \alpha \Delta u\right) ^{2}\psi
_{\lambda }^{2}dxdt\geq \frac{C}{\lambda }\dint\limits_{Q_{T}}u_{t}^{2}\psi
_{\lambda }^{2}dxdt+\frac{C}{\lambda }\dsum\limits_{i,j=1}^{n}\dint%
\limits_{Q_{T}}u_{x_{i}x_{j}}^{2}\psi _{\lambda }^{2}dxdt+
\end{equation*}%
\begin{equation}
+C\dint\limits_{Q_{T}}\left( \lambda \left( \nabla u\right) ^{2}+\lambda
^{3}u^{2}\right) \psi _{\lambda }^{2}dxdt-  \label{3.4}
\end{equation}%
\begin{equation*}
-C\left( \left\Vert \partial _{\nu }u\right\Vert _{H^{1,0}\left(
S_{T}\right) }^{2}+\left\Vert u\right\Vert _{H^{2,1}\left( S_{T}\right)
}^{2}\right) \lambda ^{2}e^{2\lambda b^{2}}-
\end{equation*}%
\begin{equation*}
-C\left( \left\Vert u\left( x,T\right) \right\Vert _{H^{1}\left( \Omega
\right) }^{2}+\left\Vert u\left( x,0\right) \right\Vert _{H^{1}\left( \Omega
\right) }^{2}\right) \lambda ^{2}\exp \left( -2\lambda \left( c^{2}\frac{%
T^{2}}{4}-b^{2}\right) \right) ,\text{ }\forall \lambda \geq \lambda _{0}.
\end{equation*}

\textbf{Proof}. We prove this theorem only for the operator $\partial
_{t}-\alpha \Delta $ since the proof for the operator $\partial _{t}+\alpha
\Delta $ is almost identical. In this proof $C=C\left( \alpha
,c,Q_{T}\right) >0$ and $\lambda _{0k}=\lambda _{0k}\left( \alpha ,c,\Omega
,T\right) ,k=1,2,3,4$ denote different numbers depending only on listed
parameters. Also, $\lambda _{04}\geq \lambda _{03}\geq \lambda _{02}\geq
\lambda _{01}\geq 1$ and we set in the end $\lambda _{0}=\lambda _{04}.$ We
take in this proof $u\in C^{3}\left( \overline{Q}_{T}\right) .$\ The case $%
u\in H^{3}\left( Q_{T}\right) $ is obtained from this one via density
arguments.

Change variables $u\Leftrightarrow v$ and express appropriate derivatives of 
$u$ via derivatives of $v$,%
\begin{equation}
\left. 
\begin{array}{c}
v=u\psi _{\lambda }=u\exp \left( \lambda \left( x_{1}^{2}-c^{2}\left(
t-T/2\right) ^{2}\right) \right) , \\ 
u=v\psi _{\lambda }^{-1}=v\exp \left( -\lambda \left( x_{1}^{2}-c^{2}\left(
t-T/2\right) ^{2}\right) \right) , \\ 
u_{t}=\left( v_{t}+2\lambda c^{2}\left( t-T/2\right) \right) \psi _{\lambda
}^{-1}, \\ 
u_{x_{1}}=\left( v_{x_{1}}-2\lambda x_{1}v\right) \psi _{\lambda }^{-1}, \\ 
u_{x_{1}x_{1}}=\left( v_{x_{1}x_{1}}-4\lambda x_{1}v_{x_{1}}\right) \psi
_{\lambda }^{-1}+ \\ 
+4\lambda ^{2}x_{1}^{2}\left( 1-1/\left( 2\lambda x_{1}\right) \right) v\psi
_{\lambda }^{-1}, \\ 
u_{x_{i}x_{i}}=v_{x_{i}x_{i}}\exp \left( -\lambda \left(
x_{1}^{2}-c^{2}\left( t-T/2\right) ^{2}\right) \right) ,i=2,...,n.%
\end{array}%
\right.  \label{3.5}
\end{equation}%
Hence, by (\ref{3.1}) and (\ref{3.5}) 
\begin{equation*}
\left( u_{t}-\alpha \Delta u\right) ^{2}\psi _{\lambda }\geq
\end{equation*}%
\begin{equation}
\geq \left( 2v_{t}+8\lambda x_{1}\alpha v_{x_{1}}\right) \left[ -\alpha
\Delta v-4\lambda ^{2}x_{1}^{2}\alpha \left( 1-1/\left( 2\lambda
x_{1}\right) -\frac{c^{2}\left( t-T/2\right) }{2\lambda x_{1}}\right) v%
\right] .  \label{3.6}
\end{equation}

\subsection{Step 1. Evaluate terms with $v_{t}$ in (\protect\ref{3.6})}

\label{sec:2.1}

We have:%
\begin{equation*}
2v_{t}\left[ -\alpha \Delta v-\left( 4\lambda ^{2}x_{1}^{2}\alpha \left(
1-1/\left( 2\lambda x_{1}\right) -\frac{c^{2}\left( t-T/2\right) }{2\lambda
x_{1}}\right) v\right) \right] =
\end{equation*}%
\begin{equation*}
=\dsum\limits_{i=1}^{n}\left( -2\alpha v_{x_{i}}v_{t}\right)
_{x_{i}}+\dsum\limits_{i=1}^{n}\left( 2\alpha v_{x_{i}}v_{tx_{i}}\right) +
\end{equation*}%
\begin{equation*}
+\left( -4\lambda ^{2}x_{1}^{2}\alpha \left( 1-1/\left( 2\lambda
x_{1}\right) -\frac{c^{2}\left( t-T/2\right) }{2\lambda x_{1}}\right)
v^{2}\right) _{t}-4c^{2}\lambda x_{1}\alpha v^{2}=
\end{equation*}%
\begin{equation*}
=-4c^{2}\lambda x_{1}\alpha v^{2}+
\end{equation*}%
\begin{equation*}
+\dsum\limits_{i=1}^{n}\left( -2\alpha v_{x_{i}}v_{t}\right) _{x_{i}}+\left(
\alpha \left( \nabla v\right) ^{2}-4\lambda ^{2}x_{1}^{2}\alpha \left(
1-1/\left( 2\lambda x_{1}\right) -\frac{c^{2}\left( t-T/2\right) }{2\lambda
x_{1}}\right) v^{2}\right) _{t}.
\end{equation*}%
Thus, coming back from $v$ to $u$ via the formula in the first line of (\ref%
{3.5}), we obtain%
\begin{equation*}
2v_{t}\left[ \alpha \Delta v-\left( 4\lambda ^{2}x_{1}^{2}\alpha \left(
1-1/\left( 2\lambda x_{1}\right) -\frac{c^{2}\left( t-T/2\right) }{2\lambda
x_{1}}\right) v\right) \right] =
\end{equation*}%
\begin{equation*}
=-4c^{2}\lambda x_{1}\alpha u^{2}\psi _{\lambda }^{2}+
\end{equation*}%
\begin{equation}
+\left( -2\alpha \left( u_{x_{1}}-2\lambda x_{1}u\right) \left(
u_{t}+2\lambda c^{2}\left( t-T/2\right) u\right) \psi _{\lambda }^{2}\right)
_{x_{1}}+  \label{3.7}
\end{equation}%
\begin{equation*}
+\dsum\limits_{i=2}^{n}\left( -2\alpha u_{x_{i}}\left( u_{t}+2\lambda
c^{2}\left( t-T/2\right) u\right) \psi _{\lambda }^{2}\right) _{x_{i}}+
\end{equation*}%
\begin{equation*}
+\left( \alpha \left( \nabla u\right) ^{2}\psi _{\lambda }-4\lambda
^{2}x_{1}^{2}\alpha \left( 1-1/\left( 2\lambda x_{1}\right) -\frac{%
c^{2}\left( t-T/2\right) }{2\lambda x_{1}}\right) u^{2}\psi _{\lambda
}^{2}\right) _{t}.
\end{equation*}

\subsection{Step 2. Evaluate terms with $v_{x_{1}}$ in (\protect\ref{3.6})}

\label{sec:2.2}

We have:%
\begin{equation*}
8\lambda x_{1}\alpha v_{x_{1}}\left[ -\alpha \Delta v-4\lambda
^{2}x_{1}^{2}\alpha \left( 1-\frac{1+c^{2}\left( t-T/2\right) }{2\lambda
x_{1}}\right) v\right] =
\end{equation*}%
\begin{equation*}
=\left( -4\lambda x_{1}\alpha ^{2}v_{x_{1}}^{2}\right) _{x_{1}}+4\lambda
\alpha ^{2}v_{x_{1}}^{2}+\dsum\limits_{i=2}^{n}\left( -8\lambda x_{1}\alpha
^{2}v_{x_{1}}v_{x_{i}}\right) _{x_{i}}+
\end{equation*}%
\begin{equation}
+\dsum\limits_{i=2}^{n}\left( 8\lambda x_{1}\alpha
^{2}v_{x_{1}x_{i}}v_{x_{i}}\right) +  \label{3.8}
\end{equation}%
\begin{equation*}
\left( -16\lambda ^{2}x_{1}^{3}\alpha ^{2}\left( 1-\frac{1+c^{2}\left(
t-T/2\right) }{2\lambda x_{1}}\right) v^{2}\right) _{x_{1}}+
\end{equation*}%
\begin{equation*}
+48\lambda ^{3}x_{1}^{2}\alpha ^{2}\left( 1-\frac{2\left( 1+c^{2}\left(
t-T/2\right) \right) }{3\lambda x_{1}}\right) v^{2}.
\end{equation*}%
Since $8\lambda x_{1}\alpha ^{2}v_{x_{1}x_{i}}v_{x_{i}}=\left( 4\lambda
x_{1}\alpha ^{2}v_{x_{i}}^{2}\right) _{x_{1}}-4\lambda \alpha
^{2}v_{x_{i}}^{2},$ then it follows from (\ref{3.8}), (\ref{3.1}) and the
first line of (\ref{3.5}) that we can choose $\lambda _{01}\geq 1$ so large
that 
\begin{equation*}
8\lambda x_{1}\alpha v_{x_{1}}\left[ -\alpha \Delta v-4\lambda
^{2}x_{1}^{2}\alpha \left( 1-\frac{1+c^{2}\left( t-T/2\right) }{2\lambda
x_{1}}\right) v\right] \geq
\end{equation*}%
\begin{equation*}
\geq -4\lambda \alpha ^{2}\left( \nabla u\right) ^{2}\psi _{\lambda
}^{2}+47\lambda ^{3}x_{1}^{2}\alpha ^{2}u^{2}\psi _{\lambda }^{2}+
\end{equation*}%
\begin{equation}
+\left( -4\lambda x_{1}\alpha ^{2}\left( u_{x_{1}}-2\lambda x_{1}u\right)
^{2}\psi _{\lambda }^{2}+\dsum\limits_{i=2}^{n}\left( 4\lambda x_{1}\alpha
^{2}u_{x_{i}}^{2}\right) \psi _{\lambda }^{2}\right) _{x_{1}}+  \label{3.9}
\end{equation}%
\begin{equation*}
+\dsum\limits_{i=2}^{n}\left( -8\lambda x_{1}\alpha ^{2}\left(
u_{x_{1}}-2\lambda x_{1}u\right) u_{x_{i}}\psi _{\lambda }^{2}\right)
_{x_{i}},\text{ }\forall \lambda \geq \lambda _{01}.
\end{equation*}%
Thus, summing up (\ref{3.7}) and (\ref{3.9}), we obtain that (\ref{3.6})
becomes with a different $\lambda _{02}\geq \lambda _{01}$: 
\begin{equation*}
\left( u_{t}-\alpha \Delta u\right) ^{2}\psi _{\lambda }^{2}\geq -4\lambda
\alpha ^{2}\left( \nabla u\right) ^{2}\psi _{\lambda }^{2}+46\lambda
^{3}x_{1}^{2}\alpha ^{2}u^{2}\psi _{\lambda }^{2}+
\end{equation*}%
\begin{equation*}
+\left( -2\alpha \left( u_{x_{1}}-2\lambda x_{1}u\right) \left(
u_{t}+2\lambda c^{2}\left( t-T/2\right) u\right) \psi _{\lambda }^{2}\right)
_{x_{1}}+
\end{equation*}%
\begin{equation*}
+\left( -4\lambda x_{1}\alpha ^{2}\left( u_{x_{1}}-2\lambda x_{1}u\right)
^{2}\psi _{\lambda }^{2}+\dsum\limits_{i=2}^{n}\left( 4\lambda x_{1}\alpha
^{2}u_{x_{i}}^{2}\right) \psi _{\lambda }^{2}\right) _{x_{1}}+
\end{equation*}%
\begin{equation}
+\dsum\limits_{i=2}^{n}\left( -2\alpha u_{x_{i}}\left( u_{t}+2\lambda
c^{2}\left( t-T/2\right) u\right) \psi _{\lambda }^{2}\right) _{x_{i}}+
\label{3.10}
\end{equation}%
\begin{equation*}
+\dsum\limits_{i=2}^{n}\left( -8\lambda x_{1}\alpha ^{2}\left(
u_{x_{1}}-2\lambda x_{1}u\right) u_{x_{i}}\psi _{\lambda }^{2}\right)
_{x_{i}}+\text{ }
\end{equation*}%
\begin{equation*}
+\left( \alpha \left( \nabla u\right) ^{2}\psi _{\lambda }^{2}-4\lambda
^{2}x_{1}^{2}\alpha \left( 1-1/\left( 2\lambda x_{1}\right) -\frac{%
c^{2}\left( t-T/2\right) }{2\lambda x_{1}}\right) u^{2}\psi _{\lambda
}^{2}\right) _{t},\text{ }\forall \lambda \geq \lambda _{02}.
\end{equation*}%
An inconvenient feature of estimate (\ref{3.10}) is the presence of the term 
$-4\lambda \alpha ^{2}\left( \nabla u\right) ^{2}\psi _{\lambda }$ with the
negative sign in the first line of (\ref{3.10}). Thus, we continue.

\subsection{Step 3. Estimate $\left( u_{t}-\protect\alpha \Delta u\right) u%
\protect\psi _{\protect\lambda }$}

\label{sec:2.3}

We have:%
\begin{equation*}
\left( u_{t}-\alpha \Delta u\right) u\psi _{\lambda }^{2}=\left( \frac{u^{2}%
}{2}\psi _{\lambda }^{2}\right) _{t}+2\lambda c^{2}\left( t-T/2\right)
u^{2}\psi _{\lambda }^{2}+\dsum\limits_{i=1}^{n}\left( -\alpha
u_{x_{i}}u\psi _{\lambda }^{2}\right) _{x_{i}}+
\end{equation*}%
\begin{equation*}
+\alpha \left( \nabla u\right) ^{2}\psi _{\lambda }^{2}+4\lambda x_{1}\alpha
u_{x_{1}}u\psi _{\lambda }^{2}=
\end{equation*}%
\begin{equation*}
=\alpha \left( \nabla u\right) ^{2}\psi _{\lambda }^{2}+\left( 2\lambda
x_{1}\alpha u^{2}\psi _{\lambda }^{2}\right) _{x_{1}}-8\lambda
^{2}x_{1}^{2}\alpha \left( 1+\frac{1}{4\lambda x_{1}^{2}}+\frac{c^{2}\left(
t-T/2\right) }{4\lambda \alpha x_{1}^{2}}\right) u^{2}\psi _{\lambda }^{2}+
\end{equation*}%
\begin{equation*}
+\dsum\limits_{i=1}^{n}\left( -\alpha u_{x_{i}}u\psi _{\lambda }^{2}\right)
_{x_{i}}.
\end{equation*}%
Thus, with another $\lambda _{03}\geq \lambda _{02}$%
\begin{equation*}
\left( u_{t}-\alpha \Delta u\right) u\psi _{\lambda }^{2}\geq \alpha \left(
\nabla u\right) ^{2}\psi _{\lambda }^{2}-9\lambda ^{2}x_{1}^{2}u^{2}\psi
_{\lambda }^{2}+
\end{equation*}%
\begin{equation}
+\dsum\limits_{i=1}^{n}\left( -\alpha u_{x_{i}}u\psi _{\lambda }^{2}\right)
_{x_{i}}+\left( 2\lambda x_{1}u^{2}\psi _{\lambda }^{2}\right)
_{x_{1}}+\left( \frac{u^{2}}{2}\psi _{\lambda }^{2}\right) _{t},\text{ }%
\forall \lambda \geq \lambda _{03}.  \label{3.11}
\end{equation}

\subsection{Step 4. Multiply (\protect\ref{3.11}) by $5\protect\lambda 
\protect\alpha $, sum up with (\protect\ref{3.10}) and integrate the
resulting inequality over $Q_{T}$}

\label{sec:2.4}

Using Gauss formula, (\ref{2.50}), (\ref{3.1}), (\ref{3.100}), (\ref{3.10}),
(\ref{3.11}) and Cauchy-Schwarz inequality, we obtain%
\begin{equation*}
5\lambda \dint\limits_{Q_{T}}\left( u_{t}-\alpha \Delta u\right) u\psi
_{\lambda }^{2}dxdt+\dint\limits_{Q_{T}}\left( u_{t}-\alpha \Delta u\right)
^{2}\psi _{\lambda }^{2}dxdt\geq
\end{equation*}%
\begin{equation}
\geq \lambda \alpha ^{2}\dint\limits_{Q_{T}}\left( \nabla u\right) ^{2}\psi
_{\lambda }^{2}dxdt+\lambda ^{3}\dint\limits_{Q_{T}}x_{1}^{2}\alpha
^{2}u^{2}\psi _{\lambda }^{2}dxdt-  \label{3.12}
\end{equation}%
\begin{equation*}
-C\left( \left\Vert \partial _{\nu }u\right\Vert _{L_{2}\left( S_{T}\right)
}^{2}+\left\Vert u\right\Vert _{H^{1,1}\left( S_{T}\right) }^{2}\right)
\lambda ^{2}e^{2\lambda b^{2}}-
\end{equation*}%
\begin{equation*}
-C\left( \left\Vert u\left( x,T\right) \right\Vert _{H^{1}\left( \Omega
\right) }^{2}+\left\Vert u\left( x,0\right) \right\Vert _{H^{1}\left( \Omega
\right) }^{2}\right) \lambda ^{2}\exp \left( -2\lambda c^{2}\left( \frac{%
T^{2}}{4}-b^{2}\right) \right) ,
\end{equation*}%
for all $\lambda \geq \lambda _{04}$ with a number $\lambda _{04}\geq
\lambda _{03}.$ By Cauchy-Schwarz inequality%
\begin{equation*}
5\lambda \left( u_{t}-\alpha \Delta u\right) u\psi _{\lambda }^{2}\leq \frac{%
5}{2}\lambda ^{2}u^{2}\psi _{\lambda }^{2}+\frac{5}{2}\left( u_{t}-\alpha
\Delta u\right) ^{2}\psi _{\lambda }^{2}.
\end{equation*}%
Hence, using (\ref{3.12}), we obtain$:$ 
\begin{equation*}
\dint\limits_{Q_{T}}\left( u_{t}-\alpha \Delta u\right) ^{2}\psi _{\lambda
}dxdt\geq C\dint\limits_{Q_{T}}\left( \lambda \left( \nabla u\right)
^{2}+\lambda ^{3}u^{2}\right) \psi _{\lambda }dxdt-
\end{equation*}%
\begin{equation}
-C\left( \left\Vert \partial _{\nu }u\right\Vert _{L_{2}\left( S_{T}\right)
}^{2}+\left\Vert u\right\Vert _{H^{1}\left( S_{T}\right) }^{2}\right)
\lambda ^{2}e^{2\lambda b^{2}}-  \label{3.13}
\end{equation}%
\begin{equation*}
-C\left( \left\Vert u\left( x,T\right) \right\Vert _{H^{1}\left( \Omega
\right) }^{2}+\left\Vert u\left( x,0\right) \right\Vert _{H^{1}\left( \Omega
\right) }^{2}\right) \lambda ^{2}\exp \left( -2\lambda \left( c^{2}\frac{%
T^{2}}{4}-b^{2}\right) \right) ,\text{ }\forall \lambda \geq \lambda _{04}.
\end{equation*}

Thus, we have obtained Carleman estimate (\ref{3.13}), which estimates lower
order derivatives $\left( \nabla u\right) ^{2}$ and $u^{2}$ in its right
hand side. We now need to incorporate derivatives which are involved in the
principal part of the operator $\partial _{t}-\alpha \Delta .$

\subsection{Step 5. Incorporating $u_{x_{i}x_{j}}^{2}$ and $u_{t}^{2}$ in
the right hand side of (\protect\ref{3.13})}

\label{sec:2.5}

We have:%
\begin{equation}
\left. 
\begin{array}{c}
\left( u_{t}-\alpha \Delta u\right) ^{2}\psi _{\lambda }^{2}=\left(
u_{t}^{2}-2\alpha u_{t}\Delta u+\alpha ^{2}\left( \Delta u\right)
^{2}\right) \psi _{\lambda }^{2}= \\ 
=u_{t}^{2}+\dsum\limits_{i=1}^{n}\left( -2\alpha u_{t}u_{x_{i}}\psi
_{\lambda }^{2}\right) _{x_{i}}+\dsum\limits_{i=1}^{n}\left( 2\alpha
u_{x_{i}t}u_{x_{i}}\psi _{\lambda }^{2}\right) +8\lambda x_{1}\alpha
u_{t}u_{x_{1}}\psi _{\lambda }^{2} \\ 
+\alpha ^{2}\left( \Delta u\right) ^{2}\psi _{\lambda }^{2}= \\ 
=u_{t}^{2}\psi _{\lambda }^{2}+\dsum\limits_{i=1}^{n}\left( -2\alpha
u_{t}u_{x_{i}}\psi _{\lambda }^{2}\right) _{x_{i}}+\left( \alpha \left(
\nabla u\right) ^{2}\psi _{\lambda }^{2}\right) _{t}+ \\ 
+4\lambda \alpha c^{2}\left( t-T/2\right) \left( \nabla u\right) ^{2}\psi
_{\lambda }^{2}+ \\ 
+8\lambda x_{1}\alpha u_{t}u_{x_{1}}\psi _{\lambda }^{2}+\alpha ^{2}\left(
\Delta u\right) ^{2}\psi _{\lambda }^{2}\geq \\ 
\geq \frac{1}{2}u_{t}^{2}\psi _{\lambda }^{2}-D\lambda ^{2}\left( \nabla
u\right) ^{2}\psi _{\lambda }^{2}+\alpha ^{2}\left( \Delta u\right) ^{2}\psi
_{\lambda }^{2}+ \\ 
+\dsum\limits_{i=1}^{n}\left( -2\alpha u_{t}u_{x_{i}}\psi _{\lambda
}^{2}\right) _{x_{i}}+\left( \alpha \left( \nabla u\right) ^{2}\psi
_{\lambda }^{2}\right) _{t}.%
\end{array}%
\right.  \label{3.130}
\end{equation}%
Here and below $D=D\left( \alpha ,c,\Omega ,T\right) >0$ denotes different
positive constants depending only on listed parameters. Thus, it follows
from (\ref{3.130}) that 
\begin{equation}
\left( u_{t}-\alpha \Delta u\right) ^{2}\psi _{\lambda }\geq \frac{1}{2}%
u_{t}^{2}\psi _{\lambda }-C_{1}\lambda ^{2}\left( \nabla u\right) ^{2}\psi
_{\lambda }+\alpha ^{2}\left( \Delta u\right) ^{2}\psi _{\lambda }+
\label{3.14}
\end{equation}%
\begin{equation*}
+\dsum\limits_{i=1}^{n}\left( -2\alpha u_{t}u_{x_{i}}\psi _{\lambda }\right)
_{x_{i}}+\left( \alpha \left( \nabla u\right) ^{2}\psi _{\lambda }\right)
_{t},\text{ }\forall \lambda \geq \lambda _{04}.
\end{equation*}

We now estimate the term $\alpha ^{2}\left( \Delta u\right) ^{2}\psi
_{\lambda }$ from the below. We have:%
\begin{equation}
\left. 
\begin{array}{c}
\alpha ^{2}\left( \Delta u\right) ^{2}\psi _{\lambda }^{2}= \\ 
=\alpha ^{2}\left( \dsum\limits_{i=1}^{n}u_{x_{i}x_{i}}^{2}\right) \psi
_{\lambda }^{2}+\alpha
^{2}\dsum\limits_{i=2}^{n}u_{x_{i}x_{i}}u_{x_{1}x_{1}}\psi _{\lambda
}^{2}+\alpha ^{2}\dsum\limits_{i,j=2,i\neq
j}^{n}u_{x_{i}x_{i}}u_{x_{j}x_{j}}\psi _{\lambda }^{2}.%
\end{array}%
\right.  \label{3.15}
\end{equation}%
Consider the second term in the second line of (\ref{3.15}), 
\begin{equation*}
\alpha ^{2}\dsum\limits_{i=2}^{n}u_{x_{i}x_{i}}u_{x_{1}x_{1}}\psi _{\lambda
}^{2}=\alpha ^{2}\dsum\limits_{i=2}^{n}\left( u_{x_{i}}u_{x_{1}x_{1}}\psi
_{\lambda }^{2}\right) _{x_{i}}-\alpha ^{2}\dsum\limits_{i=2}^{n}\left(
u_{x_{i}}u_{x_{1}x_{1}x_{i}}\psi _{\lambda }^{2}\right) =
\end{equation*}%
\begin{equation*}
=\alpha ^{2}\dsum\limits_{i=2}^{n}\left( u_{x_{i}}u_{x_{1}x_{1}}\psi
_{\lambda }^{2}\right) _{x_{i}}+\dsum\limits_{i=2}^{n}\left( -\alpha
^{2}u_{x_{i}}u_{x_{1}x_{i}}\psi _{\lambda }^{2}\right) _{x_{1}}+\alpha
^{2}\dsum\limits_{i=2}^{n}u_{x_{1}x_{i}}^{2}\psi _{\lambda }^{2}+
\end{equation*}%
\begin{equation*}
+4\alpha ^{2}\lambda x_{1}\dsum\limits_{i=2}^{n}u_{x_{i}}u_{x_{1}x_{i}}\psi
_{\lambda }^{2}\geq \frac{\alpha ^{2}}{2}\dsum%
\limits_{i=2}^{n}u_{x_{1}x_{i}}^{2}\psi _{\lambda }^{2}-C_{1}\lambda
^{2}\left( \nabla u\right) ^{2}\psi _{\lambda }^{2}+
\end{equation*}%
\begin{equation*}
+\alpha ^{2}\dsum\limits_{i=2}^{n}\left( u_{x_{i}}u_{x_{1}x_{1}}\psi
_{\lambda }^{2}\right) _{x_{i}}+\dsum\limits_{i=2}^{n}\left( -\alpha
^{2}u_{x_{i}}u_{x_{1}x_{i}}\psi _{\lambda }^{2}\right) _{x_{1}},\text{ }%
\forall \lambda \geq \lambda _{04}.
\end{equation*}%
Thus, the second term in the second line of (\ref{3.15}) can be estimated as:%
\begin{equation}
\alpha ^{2}\dsum\limits_{i=2}^{n}u_{x_{i}x_{i}}u_{x_{1}x_{1}}\psi _{\lambda
}^{2}\geq \frac{\alpha ^{2}}{2}\dsum\limits_{i=2}^{n}u_{x_{1}x_{i}}^{2}\psi
_{\lambda }^{2}-C\lambda ^{2}\left( \nabla u\right) ^{2}\psi _{\lambda }^{2}+
\label{3.16}
\end{equation}%
\begin{equation*}
+\alpha ^{2}\dsum\limits_{i=2}^{n}\left( u_{x_{i}}u_{x_{1}x_{1}}\psi
_{\lambda }^{2}\right) _{x_{i}}+\dsum\limits_{i=2}^{n}\left( -\alpha
^{2}u_{x_{i}}u_{x_{1}x_{i}}\psi _{\lambda }^{2}\right) _{x_{1}},\text{ }%
\forall \lambda \geq \lambda _{04}.
\end{equation*}

We now estimate the third term in the second line of (\ref{3.15}),%
\begin{equation*}
\left. 
\begin{array}{c}
\alpha ^{2}\dsum\limits_{i,j=2,i\neq j}^{n}u_{x_{i}x_{i}}u_{x_{j}x_{j}}\psi
_{\lambda }^{2}=\dsum\limits_{i,j=2,i\neq j}^{n}\left( \alpha
^{2}u_{x_{i}x_{i}}u_{x_{j}}\psi _{\lambda }^{2}\right)
_{x_{j}}-\dsum\limits_{i,j=2,i\neq j}^{n}\left( \alpha
^{2}u_{x_{i}x_{i}x_{j}}u_{x_{j}}\psi _{\lambda }^{2}\right) = \\ 
=\dsum\limits_{i,j=2,i\neq j}^{n}\left( \alpha
^{2}u_{x_{i}x_{i}}u_{x_{j}}\psi _{\lambda }^{2}\right)
_{x_{j}}+\dsum\limits_{i,j=2,i\neq j}^{n}\left( -\alpha
^{2}u_{x_{i}x_{j}}u_{x_{j}}\psi _{\lambda }^{2}\right) _{x_{i}}+\alpha
^{2}\dsum\limits_{i,j=2,i\neq j}^{n}u_{x_{i}x_{j}}^{2}\psi _{\lambda }^{2},%
\text{ }\forall \lambda \geq \lambda _{04}.%
\end{array}%
\right.
\end{equation*}%
Thus,%
\begin{equation}
\left. 
\begin{array}{c}
\alpha ^{2}\dsum\limits_{i,j=2,i\neq j}^{n}u_{x_{i}x_{i}}u_{x_{j}x_{j}}\psi
_{\lambda }^{2}=\alpha ^{2}\dsum\limits_{i,j=2,i\neq
j}^{n}u_{x_{i}x_{j}}^{2}\psi _{\lambda }^{2}+ \\ 
+\dsum\limits_{i,j=2,i\neq j}^{n}\left( \alpha
^{2}u_{x_{i}x_{i}}u_{x_{j}}\psi _{\lambda }^{2}\right)
_{x_{j}}+\dsum\limits_{i,j=2,i\neq j}^{n}\left( -\alpha
^{2}u_{x_{i}x_{j}}u_{x_{j}}\psi _{\lambda }^{2}\right) _{x_{i}},\text{ }%
\forall \lambda \geq \lambda _{04}.%
\end{array}%
\right.  \label{3.17}
\end{equation}%
Using (\ref{3.15})-(\ref{3.17}), we obtain%
\begin{equation*}
\alpha ^{2}\left( \Delta u\right) ^{2}\psi _{\lambda }^{2}\geq \frac{\alpha
^{2}}{2}\dsum\limits_{i=2}^{n}u_{x_{i}x_{j}}^{2}\psi _{\lambda
}^{2}-D\lambda ^{2}\left( \nabla u\right) ^{2}\psi _{\lambda }^{2}+
\end{equation*}%
\begin{equation*}
+\dsum\limits_{i,j=1}^{n}\left( \alpha ^{2}u_{x_{i}x_{i}}u_{x_{j}}\psi
_{\lambda }^{2}\right) _{x_{j}}+\dsum\limits_{i,j=1}^{n}\left( -\alpha
^{2}u_{x_{i}x_{j}}u_{x_{j}}\psi _{\lambda }^{2}\right) _{x_{i}}.
\end{equation*}%
Comparing this with (\ref{3.14}), we obtain%
\begin{equation*}
\left( u_{t}-\alpha \Delta u\right) ^{2}\psi _{\lambda }^{2}\geq \left( 
\frac{1}{2}u_{t}^{2}\psi _{\lambda }^{2}+\frac{\alpha ^{2}}{2}%
\dsum\limits_{i=2}^{n}u_{x_{i}x_{j}}^{2}\psi _{\lambda }^{2}\right)
-D\lambda ^{2}\left( \nabla u\right) ^{2}\psi _{\lambda }^{2}+
\end{equation*}%
\begin{equation}
+\dsum\limits_{i=1}^{n}\left( -2\alpha u_{t}u_{x_{i}}\psi _{\lambda
}^{2}\right) _{x_{i}}+\left( \alpha \left( \nabla u\right) ^{2}\psi
_{\lambda }^{2}\right) _{t}+  \label{3.18}
\end{equation}%
\begin{equation*}
+\alpha ^{2}\dsum\limits_{i=2}^{n}\left( u_{x_{i}}u_{x_{1}x_{1}}\psi
_{\lambda }^{2}\right) _{x_{i}}+\dsum\limits_{i=2}^{n}\left( -\alpha
^{2}u_{x_{i}}u_{x_{1}x_{i}}\psi _{\lambda }^{2}\right) _{x_{1}}
\end{equation*}%
\begin{equation*}
+\dsum\limits_{i,j=2,i\neq j}^{n}\left( \alpha
^{2}u_{x_{i}x_{i}}u_{x_{j}}\psi _{\lambda }^{2}\right)
_{x_{j}}+\dsum\limits_{i,j=2,i\neq j}^{n}\left( -\alpha
^{2}u_{x_{i}x_{j}}u_{x_{j}}\psi _{\lambda }^{2}\right) _{x_{i}},\text{ }%
\forall \lambda \geq \lambda _{04}.
\end{equation*}%
Set $\lambda _{0}=\lambda _{04}=\lambda _{0}\left( \alpha ,c,\Omega
,T\right) \geq 1.$ Let $\lambda \geq \lambda _{0}.$ Multiply inequality by (%
\ref{3.18}) by $C/\left( 2\lambda D\right) ,$ integrate over $Q_{T},$ use (%
\ref{3.100}) and sum up with (\ref{3.13}). Then we obtain the target
Carleman estimate (\ref{3.4}). $\ \square $

\section{H\"{o}lder Stability and Uniqueness}

\label{sec:5}

The parameter $T>0$ is actually the observation time of the game process
from the boundary $\partial \Omega $. Therefore, it makes sense to consider
two cases:

\begin{enumerate}
\item $T>0$ is an arbitrary number.

\item $T>0$ is sufficiently large as compared with $\sqrt{b^{2}-a^{2}},$
which is a number associated with the length $\left( b-a\right) $ of the
prism $\Omega $ in the $x_{1}-$direction, see (\ref{2.1}).
\end{enumerate}

In particular, it follows from Theorem 5.2 that the problem we consider
admits a stronger stability estimate for sufficiently large $T$ than for an
arbitrary $T$, which is to be expected, see Remarks 5.2.

\subsection{ The case of an arbitrary $T>0$}

\label{sec:5.1}

Consider an arbitrary number $\rho $ such that 
\begin{equation}
\rho \in \left( 0,\frac{\sqrt{13}-3}{2}\right) .  \label{4.01}
\end{equation}%
Elementary calculations show that: 
\begin{equation}
\left. 
\begin{array}{c}
1-\rho ^{2}-3\rho >0\text{ for }\rho \text{ satisfying (\ref{4.01}),} \\ 
1/5<\left( \sqrt{13}-3\right) /2. \\ 
\text{If }\rho \in \left( 1/5,\left( \sqrt{13}-3\right) /2\right) ,\text{
then} \\ 
T/\left( 2\left( \rho +3\right) \right) \left( 1+\sqrt{1-3\rho -\rho ^{2}}%
\right) <T/4\text{.} \\ 
\text{If }\rho \in \left( 0,1/5\right) ,\text{ then} \\ 
T/\left( 2\left( \rho +3\right) \right) \left( 1+\sqrt{1-3\rho -\rho ^{2}}%
\right) >T/4\text{ .} \\ 
\text{If }\rho \in \left( 0,\left( \sqrt{13}-3\right) /2\right) ,\text{ then}
\\ 
T/\left( 2\left( \rho +3\right) \right) \left( 1+\sqrt{1-3\rho -\rho ^{2}}%
\right) <T/4\text{.}%
\end{array}%
\right.  \label{4.02}
\end{equation}
Choose a number $\varepsilon >0.$ By (\ref{3.01}) we want $\varepsilon \in
\left( 0,T/4\right) .$ Hence, keeping in mind (\ref{4.01}) and (\ref{4.02}),
let the number $\varepsilon $ be such that: 
\begin{equation}
\left. 
\begin{array}{c}
\text{If }\rho \in \left( 1/5,\left( \sqrt{13}-3\right) /2\right) ,\text{
then} \\ 
T/\left( 2\left( \rho +3\right) \right) \left( 1-\sqrt{1-3\rho -\rho ^{2}}%
\right) <\varepsilon < \\ 
<T/\left( 2\left( \rho +3\right) \right) \left( 1+\sqrt{1-3\rho -\rho ^{2}}%
\right) <T/4. \\ 
\text{If }\rho \in \left( 0,1/5\right) ,\text{ then} \\ 
T/\left( 2\left( \rho +3\right) \right) \left( 1-\sqrt{1-3\rho -\rho ^{2}}%
\right) <\varepsilon <T/4.%
\end{array}%
\right.  \label{4.03}
\end{equation}

\textbf{Theorem 5.1 }(arbitrary $\dot{T}>0$):

\begin{enumerate}
\item \emph{Choose an arbitrary number }$\rho $\emph{\ satisfying (\ref{4.01}%
) and let }$\gamma \in \left( 0,\rho \right) $\emph{\ be an arbitrary
number. Let }$\varepsilon $\emph{\ be another arbitrary number satisfying (%
\ref{4.03}). }

\item \emph{Let the kernel }$K\left( x,y\right) $\emph{\ of the integral
operator in (\ref{2.6}) has either of the forms (\ref{30.02}), (\ref{30.03})
or \ref{30.05}), (\ref{30.06}). Let }$R_{1}>0$\emph{\ be the number defined
in (\ref{30.03}) and (\ref{30.06}). }

\item \emph{Let} $R_{2},R_{3},R_{4},R_{5}>0$\emph{\ be some numbers. Let the
function }$K\left( x,y\right) $ \emph{satisfies either conditions of
subsection 3.1 or conditions of subsection 3.2. In (\ref{2.6}), let the
function }$P=P\left( x,t,z_{1},z_{2}\right) :\overline{Q}_{T}\times \mathbb{R%
}^{2}\rightarrow \mathbb{R}$\emph{\ be bounded in any bounded subset of the
set }$\overline{Q}_{T}\times \mathbb{R}^{2}$\emph{\ and such that there
exist derivatives }$P_{z_{1}},P_{z_{2}}\in C\left( \overline{Q}_{T}\times 
\mathbb{R}^{2}\right) $ \emph{satisfying} 
\begin{equation}
\max \left( \sup_{Q_{T}\times \mathbb{R}^{2}}\left\vert P_{z_{1}}\left(
x,t,y,z\right) \right\vert ,\sup_{Q_{T}\times \mathbb{R}^{2}}\left\vert
P_{z_{2}}\left( x,t,y,z\right) \right\vert \right) \leq R_{2}.  \label{4.3}
\end{equation}

\item \emph{In (\ref{2.6}), let the function }$k\in C^{1}\left( \overline{%
\Omega }\right) $ \emph{and}%
\begin{equation}
\left\Vert k\right\Vert _{C^{1}\left( \overline{\Omega }\right) }\leq R_{3}.
\label{4.4}
\end{equation}%
\emph{Consider two sets of functions }$K_{4}\left( R_{4}\right) $ \emph{and }%
$K_{5}\left( R_{5}\right) $\emph{,}%
\begin{equation}
K_{4}\left( R_{4}\right) =\left\{ u\in H^{3}\left( Q_{T}\right)
:\sup_{Q_{T}}\left\vert u\right\vert ,\sup_{Q_{T}}\left\vert \nabla
u\right\vert ,\sup_{Q_{T}}\left\vert \Delta u\right\vert \leq R_{4}\right\} ,
\label{4.5}
\end{equation}%
\begin{equation}
K_{5}\left( R_{5}\right) =\left\{ m\in H^{3}\left( Q_{T}\right)
:\sup_{Q_{T}}\left\vert m\right\vert ,\sup_{Q_{T}}\left\vert \nabla
m\right\vert \leq R_{5}\right\} .  \label{4.6}
\end{equation}%
\emph{Introduce the number }$R$\emph{,}%
\begin{equation}
R=\max \left( R_{1},R_{2},R_{3},R_{4},R_{5}\right) .  \label{4.7}
\end{equation}%
\emph{Let two vector functions }%
\begin{equation}
\left( u_{1},m_{1}\right) ,\left( u_{2},m_{2}\right) \in K_{4}\left(
R_{4}\right) \times K_{5}\left( R_{5}\right)  \label{4.8}
\end{equation}%
\emph{\ satisfy equations (\ref{2.6}) as well as the following Dirichlet and
Neumann boundary conditions (see (\ref{2.8}), (\ref{2.9}) and (\ref{2.13})):}%
\begin{equation}
u_{i}\mid _{S_{T}}=g_{0,i},\text{ }\partial _{\nu }u_{i}\mid
_{S_{T}}=g_{1,i},\text{ }i=1,2,  \label{4.9}
\end{equation}%
\begin{equation}
m_{i}\mid _{S_{T}}=p_{0,i},\text{ }\partial _{\nu }m_{i}\mid
_{S_{T}}=p_{1,i},\text{ }i=1,2.  \label{4.10}
\end{equation}

\item \emph{Let }$\delta \in \left( 0,1\right) $\emph{\ be a sufficiently
small number characterizing the level of the error in the boundary data (\ref%
{4.9}), (\ref{4.10}). Assume that}%
\begin{equation}
\left. 
\begin{array}{c}
\left\Vert g_{0,1}-g_{0,2}\right\Vert _{H^{2,1}\left( S_{T}\right)
},\left\Vert g_{1,1}-g_{0,2}\right\Vert _{H^{1,0}\left( S_{T}\right) }\leq
\delta , \\ 
\left\Vert p_{0,1}-p_{0,2}\right\Vert _{H^{2,1}\left( S_{T}\right)
},\left\Vert p_{1,1}-p_{0,2}\right\Vert _{H^{1,0}\left( S_{T}\right) }\leq
\delta .%
\end{array}%
\right.  \label{4.11}
\end{equation}%
\emph{\ \ }

\item \emph{Let }$Q_{T,\varepsilon }$\emph{\ and }$Q_{T,2\varepsilon }$\emph{%
\ be two subdomains of the domain }$Q_{T}$\emph{, which are defined in (\ref%
{3.01}). }
\end{enumerate}

\emph{Then:}

\textbf{A.} \emph{For any pair of numbers }$\rho $\emph{\ and }$\varepsilon $%
\emph{\ satisfying satisfying (\ref{4.01}) and (\ref{4.03}) respectively
there exists a sufficiently small number }$\delta _{1}=\delta _{1}\left(
\alpha ,\rho ,\gamma ,\varepsilon ,R,Q_{T}\right) \in \left( 0,1\right) $%
\emph{\ depending only on listed parameters such that the following H\"{o}%
lder stability estimate holds:}%
\begin{equation}
\left. 
\begin{array}{c}
\left\Vert u_{1}-u_{2}\right\Vert _{H^{2,1}\left( Q_{T,2\varepsilon }\right)
}+\left\Vert m_{1}-m_{2}\right\Vert _{H^{2,1}\left( Q_{T,2\varepsilon
}\right) }\leq \\ 
\leq C_{1}\left( 1+\left\Vert u_{1}-u_{2}\right\Vert _{L_{2}\left(
Q_{T}\right) }+\left\Vert m_{1}-m_{2}\right\Vert _{H^{2}\left( Q_{T}\right)
}\right) \delta ^{\rho -\gamma },\text{ }\forall \delta \in \left( 0,\delta
_{1}\right) ,\text{ }%
\end{array}%
\right.  \label{4.13}
\end{equation}%
\emph{where the number }$C_{1}=C_{1}\left( \alpha ,\rho ,\varepsilon
,R,Q_{T}\right) >0$\emph{\ depends only on listed parameters.}

\textbf{B}. \emph{If in (\ref{4.9})and (\ref{4.10}) }%
\begin{equation}
g_{0,1}\equiv g_{0,2},\text{ }g_{1,1}\equiv g_{1,2},\text{ }p_{0,1}\equiv
p_{0,2},\text{ }p_{1,1}\equiv p_{1,2},\text{ }\left( x,t\right) \in S_{T},
\label{4.140}
\end{equation}%
\emph{then }%
\begin{equation}
u_{1}\left( x,t\right) \equiv u_{2}\left( x,t\right) \text{ and }m_{1}\left(
x,t\right) \equiv m_{2}\left( x,t\right) ,\text{ }\left( x,t\right) \in
Q_{T},  \label{4.141}
\end{equation}%
\emph{which means that there exists at most one pair of functions }$\left(
u,m\right) \in K_{4}\left( R_{4}\right) \times K_{5}\left( R_{5}\right) $%
\emph{\ satisfying MFGS system (\ref{2.6}) and boundary conditions (\ref{2.9}%
).}

\textbf{Remarks 5.1:}

\begin{enumerate}
\item \emph{We remind the statement, which was made in section 1: it follows
from Theorem 5.1 that, unlike the conventional case \cite{A}, we do not need
to know for the uniqueness functions }$u\left( x,T\right) $\emph{\ and }$%
m\left( x,0\right) .$

\item \emph{It is typical in the field of Ill-Posed and Inverse Problems to
impose a priori bounds like ones in (\ref{4.8}): that solution belongs to an
a priori chosen bounded set, see, e.g. \cite{T}.}

\item \emph{To simplify the presentation, we use in the proof of this
theorem only the first form (\ref{30.02}), (\ref{30.03}) of the kernel }$%
K\left( x,y\right) .$\emph{\ The second form can be considered absolutely
similarly due to Lemma 3.1.}
\end{enumerate}

\textbf{Proof of Theorem 5.1}. We set in (\ref{3.1}) 
\begin{equation}
c^{2}=c^{2}\left( \rho ,\varepsilon ,T\right) =\frac{\left( 1+\beta \right)
\left( b^{2}-a^{2}\right) }{\varepsilon \left( T-\varepsilon \right) },\text{
}\beta >0,  \label{4.142}
\end{equation}%
where $\beta =\beta \left( \rho ,\varepsilon ,T\right) >0$ is a certain
number, which we will define later. Hence, we use in this proof $c^{2}$ from
(\ref{4.142}), which makes the dependence of the constant $C$ in estimate (%
\ref{3.4}) of Theorem 4.1 as $C=C\left( \alpha ,\rho ,\varepsilon
,Q_{T}\right) >0.$ Therefore, in this proof $C_{1}=C_{1}\left( \alpha ,\rho
,\varepsilon ,R,Q_{T}\right) >0$ and $\lambda _{1}=\lambda _{1}\left( \alpha
,\rho ,\varepsilon ,R,Q_{T}\right) $ denote different numbers depending only
on listed parameters, where $\lambda _{1}=\lambda _{1}\left( \alpha ,\rho
,\varepsilon ,R,Q_{T}\right) \geq \lambda _{0}=\lambda _{0}\left( \alpha
,\rho ,\varepsilon ,Q_{T}\right) \geq 1$. Recall that the parameter $\lambda
_{0}$ was chosen in Theorem 4.1.

Consider arbitrary numbers $d_{1},k_{1},d_{2},k_{2}\in \mathbb{R}.$ Denote $%
\widetilde{d}=d_{1}-d_{2}$ and $\widetilde{k}=k_{1}-k_{2}.$ Then 
\begin{equation}
d_{1}k_{1}-d_{2}k_{2}=\widetilde{d}k_{1}+\widetilde{k}d_{2}.  \label{4.15}
\end{equation}

Denote%
\begin{equation}
\left. 
\begin{array}{c}
\widetilde{v}\left( x,t\right) =\left( u_{1}-u_{2}\right) \left( x,t\right) ,%
\text{ }\widetilde{m}\left( x,t\right) =\left( m_{1}-m_{2}\right) \left(
x,t\right) ,\left( x,t\right) \in Q_{T}, \\ 
\widetilde{g}_{0}\left( x,t\right) =\left( g_{0,1}-g_{0,2}\right) \left(
x,t\right) ,\text{ }\widetilde{g}_{1}\left( x,t\right) =\left(
g_{1,1}-g_{1,2}\right) \left( x,t\right) ,\text{ }\left( x,t\right) \in
S_{T}, \\ 
\widetilde{p}_{0}\left( x,t\right) =\left( p_{0,1}-p_{0,2}\right) \left(
x,t\right) ,\text{ }\widetilde{p}_{1}\left( x,t\right) =\left(
p_{1,1}-p_{1,2}\right) \left( x,t\right) ,\text{ }\left( x,t\right) \in
S_{T}.%
\end{array}%
\right.  \label{4.16}
\end{equation}%
It follows from (\ref{4.11}) and (\ref{4.16}) that 
\begin{equation}
\left\Vert \widetilde{g}_{0}\right\Vert _{H^{2,1}\left( S_{T}\right)
},\left\Vert \widetilde{g}_{1}\right\Vert _{H^{1,0}\left( S_{T}\right)
},\left\Vert \widetilde{p}_{0}\right\Vert _{H^{2,1}\left( S_{T}\right)
},\left\Vert \widetilde{p}_{1}\right\Vert _{H^{1,0}\left( S_{T}\right) }\leq
\delta .  \label{4.016}
\end{equation}

The multidimensional analog of Taylor formula \cite{V} as well as (\ref{2.6}%
), item 3 of Remarks 5.1, (\ref{30.04}), (\ref{4.3}), (\ref{4.7}), (\ref%
{4.15}) and (\ref{4.16}) imply:%
\begin{equation}
\left. 
\begin{array}{c}
P\left( x,t,\dint\limits_{\Omega }K\left( x,y\right) m_{1}\left( y,t\right)
dy,m_{1}\left( x,t\right) \right) - \\ 
-P\left( x,t,\dint\limits_{\Omega }K\left( x,y\right) m_{2}\left( y,t\right)
dy,m_{2}\left( x,t\right) \right) = \\ 
=D_{1}\left( x,t\right) \dint\limits_{\Omega _{1}}K_{1}\left( x,\overline{y}%
\right) \widetilde{m}\left( x_{1},\overline{y},t\right) +D_{2}\left(
x,t\right) \widetilde{m}\left( x,t\right) ,%
\end{array}%
\right.  \label{4.17}
\end{equation}%
where functions $D_{1},D_{2}\in C\left( \overline{Q}_{T}\right) $ and 
\begin{equation}
\left\vert D_{1}\left( x,t\right) \right\vert ,\left\vert D_{2}\left(
x,t\right) \right\vert \leq R.  \label{4.18}
\end{equation}

Write equations (\ref{2.6}) for each of two pairs $\left( u_{1},m_{1}\right)
,\left( u_{2},m_{2}\right) .$ Then subtract equations for the second pair
from corresponding equations for the first pair. Then use formulas (\ref{4.4}%
)-(\ref{4.8}) and (\ref{4.15})-(\ref{4.18}). Recalling that Carleman
estimates can work with both equations and inequalities \cite{KL,LRS,Yam},
we turn two resulting equations in a more general form of inequalities,%
\begin{equation}
\left\vert \widetilde{v}_{t}+\alpha \Delta \widetilde{v}\right\vert \leq
C_{1}\left( \left\vert \nabla \widetilde{v}\right\vert +\dint\limits_{\Omega
_{1}}\left\vert \widetilde{m}\left( x_{1},\overline{y},t\right) \right\vert d%
\overline{y}+\left\vert \widetilde{m}\right\vert \right) ,\text{ }\left(
x,t\right) \in Q_{T},  \label{4.19}
\end{equation}%
\begin{equation}
\left\vert \widetilde{m}_{t}-\alpha \Delta \widetilde{m}\right\vert \leq
C_{1}\left( \left\vert \nabla \widetilde{m}\right\vert +\left\vert 
\widetilde{m}\right\vert +\left\vert \Delta \widetilde{v}\right\vert
+\left\vert \nabla \widetilde{v}\right\vert \right) ,\text{ }\left(
x,t\right) \in Q_{T},  \label{4.20}
\end{equation}%
\begin{equation}
\widetilde{v}\mid _{S_{T}}=\widetilde{g}_{0},\text{ }\partial _{\nu }%
\widetilde{v}\mid _{S_{T}}=\widetilde{g}_{1},  \label{4.21}
\end{equation}%
\begin{equation}
\widetilde{m}\mid _{S_{T}}=\widetilde{p}_{0},\text{ }\partial _{\nu }%
\widetilde{m}\mid _{S_{T}}=\widetilde{p}_{1}.  \label{4.22}
\end{equation}%
Consider now a cut-off function $\chi _{\varepsilon }\left( t\right) \in
C^{1}\left[ 0,T\right] $ satisfying the following conditions:%
\begin{equation}
\chi _{\varepsilon }\left( t\right) =\left\{ 
\begin{array}{c}
1,t\in \left( 2\varepsilon ,T-2\varepsilon \right) , \\ 
0,t\in \left( 0,\varepsilon \right) \cup \left( T-\varepsilon ,T\right) , \\ 
\text{between 0 and 1 for }t\in \left( \varepsilon ,2\varepsilon \right)
\cup \left( T-2\varepsilon ,T-\varepsilon \right) .%
\end{array}%
\right.  \label{4.220}
\end{equation}%
The existence of such functions is well known from the standard Mathematical
Analysis course. Introduce new functions $w\left( x,t\right) $ and $q\left(
x,t\right) ,$%
\begin{equation}
w\left( x,t\right) =\chi _{\varepsilon }\left( t\right) \widetilde{v}\left(
x,t\right) ,\text{ }q\left( x,t\right) =\chi _{\varepsilon }\left( t\right) 
\widetilde{m}\left( x,t\right) .  \label{4.221}
\end{equation}%
It follows from (\ref{3.01}), (\ref{4.220}) and (\ref{4.221}) that%
\begin{equation}
\left. 
\begin{array}{c}
w\left( x,t\right) =\widetilde{v}\left( x,t\right) ,\text{ }q\left(
x,t\right) =\widetilde{m}\left( x,t\right) \text{ in }Q_{\varepsilon T}, \\ 
\chi _{\varepsilon }\left( t\right) \widetilde{v}_{t}=w_{t}-\chi
_{\varepsilon }^{\prime }\left( t\right) \widetilde{v},\text{ }\chi
_{\varepsilon }\left( t\right) \widetilde{m}_{t}=q_{t}-\chi _{\varepsilon
}^{\prime }\left( t\right) \widetilde{m}, \\ 
\chi _{\varepsilon }\left( t\right) \nabla \widetilde{v}=\nabla w,\text{ }%
\chi _{\varepsilon }\left( t\right) \Delta \widetilde{v}=\Delta w, \\ 
\chi _{\varepsilon }\left( t\right) \nabla \widetilde{m}=\nabla q,\text{ }%
\chi _{\varepsilon }\left( t\right) \Delta \widetilde{m}=\Delta q.%
\end{array}%
\right.  \label{4.222}
\end{equation}%
Hence, multiplying equations (\ref{4.19})-(\ref{4.22}) by $\chi
_{\varepsilon }\left( t\right) $ and using (\ref{4.220})-(\ref{4.222}), we
obtain%
\begin{equation}
\left\vert w_{t}+\alpha \Delta w\right\vert \leq C_{1}\left( \left\vert
\nabla w\right\vert +\dint\limits_{\Omega _{1}}\left\vert q\left( x_{1},%
\overline{y},t\right) \right\vert d\overline{y}+\left\vert q\right\vert
\right) +\left\vert \chi _{\varepsilon }^{\prime }\left( t\right)
\right\vert \left\vert \widetilde{v}\right\vert ,\text{ }\left( x,t\right)
\in Q_{T},  \label{4.223}
\end{equation}%
\begin{equation}
\left\vert q_{t}-\alpha \Delta q\right\vert \leq C_{1}\left( \left\vert
\nabla q\right\vert +\left\vert q\right\vert +\left\vert \Delta w\right\vert
+\left\vert \nabla w\right\vert \right) +\left\vert \chi _{\varepsilon
}^{\prime }\left( t\right) \right\vert \left\vert \widetilde{v}\right\vert ,%
\text{ }\left( x,t\right) \in Q_{T},  \label{4.224}
\end{equation}%
\begin{equation}
w\mid _{S_{T}}=\chi _{\varepsilon }\left( t\right) \widetilde{g}_{0},\text{ }%
\partial _{\nu }w\mid _{S_{T}}=\chi _{\varepsilon }\left( t\right) 
\widetilde{g}_{1},  \label{4.225}
\end{equation}%
\begin{equation}
q\mid _{S_{T}}=\chi _{\varepsilon }\left( t\right) \widetilde{p}_{0},\text{ }%
\partial _{\nu }q\mid _{S_{T}}=\chi _{\varepsilon }\left( t\right) 
\widetilde{p}_{1}.  \label{4.226}
\end{equation}

Square both sides of each equation (\ref{4.223}) and (\ref{4.224}), apply
Cauchy-Schwarz inequality, multiply by the CWF $\psi _{\lambda }\left(
x_{1},t\right) $ defined in (\ref{3.1}) and integrate over $Q_{T}.$ Note
that since the function $\psi _{\lambda }\left( x_{1},t\right) $ is
independent on $\overline{x},$ then%
\begin{equation}
\dint\limits_{Q_{T}}\left( \dint\limits_{\Omega _{1}}\left\vert q\left(
x_{1},\overline{y},t\right) \right\vert d\overline{y}\right) ^{2}\psi
_{\lambda }dxdt\leq C_{1}\dint\limits_{Q_{T}}q^{2}\psi _{\lambda }dxdt.
\label{4.022}
\end{equation}

Thus, we obtain two integral inequalities,%
\begin{equation}
\left. 
\begin{array}{c}
\dint\limits_{Q_{T}}\left( w_{t}+\alpha \Delta w\right) ^{2}\psi _{\lambda
}dxdt\leq C_{1}\dint\limits_{Q_{T}}\left( \left( \nabla w\right)
^{2}+q^{2}\right) \psi _{\lambda }dxdt+ \\ 
+2\dint\limits_{Q_{T}}\left\vert \chi _{\varepsilon }^{\prime }\left(
t\right) \right\vert ^{2}\widetilde{v}^{2}\psi _{\lambda }dxdt,%
\end{array}%
\right.  \label{4.23}
\end{equation}%
\begin{equation}
\left. 
\begin{array}{c}
\dint\limits_{Q_{T}}\left( q_{t}-\alpha \Delta q\right) ^{2}\psi _{\lambda
}dxdt\leq C_{1}\dint\limits_{Q_{T}}\left( \left( \nabla q\right)
^{2}+q^{2}+\left( \Delta w\right) ^{2}+w^{2}\right) \psi _{\lambda }dxdt+ \\ 
+2\dint\limits_{Q_{T}}\left\vert \chi _{\varepsilon }^{\prime }\left(
t\right) \right\vert ^{2}\widetilde{m}^{2}\psi _{\lambda }dxdt.%
\end{array}%
\right.  \label{4.24}
\end{equation}

Apply Carleman estimate (\ref{3.4}) of Theorem 4.1 for the operator $%
\partial _{t}-\alpha \Delta $ to the left hand side of (\ref{4.24}) and use (%
\ref{4.11}), (\ref{4.16}), (\ref{4.016}), (\ref{4.225}) and (\ref{4.226}).
Note that (\ref{4.220}) and (\ref{4.221}) imply that the last line of (\ref%
{3.4}) equals zero if the function $u$ is replaced with the function $q$. We
obtain%
\begin{equation}
\left. 
\begin{array}{c}
\left( 1/\lambda \right) \dint\limits_{Q_{T}}q_{t}^{2}\psi _{\lambda
}dxdt+\left( 1/\lambda \right)
\dsum\limits_{i,j=1}^{n}\dint\limits_{Q_{T}}q_{x_{i}x_{j}}^{2}\psi _{\lambda
}dxdt+ \\ 
+\dint\limits_{Q_{T}}\left( \lambda \left( \nabla q\right) ^{2}+\lambda
^{3}q^{2}\right) \psi _{\lambda }dxdt\leq \\ 
\leq C_{1}\dint\limits_{Q_{T}}\left[ \left( \nabla q\right)
^{2}+q^{2}+\left( \Delta w\right) ^{2}+\left( \nabla w\right) ^{2}\right]
\psi _{\lambda }dxdt+ \\ 
+2\dint\limits_{Q_{T}}\left\vert \chi _{\varepsilon }^{\prime }\left(
t\right) \right\vert ^{2}\widetilde{m}^{2}\psi _{\lambda }dxdt+C_{1}\delta
^{2}\lambda ^{2}e^{2\lambda b^{2}},\text{ }\forall \lambda \geq \lambda _{0}.%
\end{array}%
\right.  \label{4.25}
\end{equation}%
Let $\lambda _{1}$ be so large that $\lambda _{1}\geq 2C_{1}.$ Then (\ref%
{4.25}) leads to:%
\begin{equation}
\left. 
\begin{array}{c}
\left( 1/\lambda \right) \dint\limits_{Q_{T}}q_{t}^{2}\psi _{\lambda
}dxdt+\left( 1/\lambda \right)
\dsum\limits_{i,j=1}^{n}\dint\limits_{Q_{T}}q_{x_{i}x_{j}}^{2}\psi _{\lambda
}dxdt+ \\ 
+\dint\limits_{Q_{T}}\left[ \lambda \left( \nabla q\right) ^{2}+\lambda
^{3}q^{2}\right] \psi _{\lambda }dxdt\leq \\ 
\leq C_{1}\dint\limits_{Q_{T}}\left[ \left( \Delta w\right) ^{2}+\left(
\nabla w\right) ^{2}\right] \psi _{\lambda
}dxdt+2\dint\limits_{Q_{T}}\left\vert \chi _{\varepsilon }^{\prime }\left(
t\right) \right\vert ^{2}\widetilde{m}^{2}\psi _{\lambda }dxdt+ \\ 
+C_{1}\delta ^{2}\lambda ^{2}e^{2\lambda b^{2}},\text{ }\forall \lambda \geq
\lambda _{1}.%
\end{array}%
\right.  \label{4.26}
\end{equation}%
In particular, using (\ref{4.26}), we obtain%
\begin{equation}
\left. 
\begin{array}{c}
\dint\limits_{Q_{T}}q^{2}\psi _{\lambda }dxdt\leq C_{1}\left( 1/\lambda
^{3}\right) \dint\limits_{Q_{T}}\left[ \left( \Delta w\right) ^{2}+\left(
\nabla w\right) ^{2}\right] \psi _{\lambda }dxdt+ \\ 
+2\left( 1/\lambda ^{3}\right) \dint\limits_{Q_{T}}\left\vert \chi
_{\varepsilon }^{\prime }\left( t\right) \right\vert ^{2}\widetilde{m}%
^{2}\psi _{\lambda }dxdt+C_{1}\delta ^{2}\lambda ^{2}e^{2\lambda b^{2}},%
\text{ }\forall \lambda \geq \lambda _{1}.%
\end{array}%
\right.  \label{4.27}
\end{equation}%
Substituting (\ref{4.27}) in (\ref{4.23}), we obtain%
\begin{equation}
\left. 
\begin{array}{c}
\dint\limits_{Q_{T}}\left( w_{t}+\alpha \Delta w\right) ^{2}\psi _{\lambda
}dxdt\leq C_{1}\dint\limits_{Q_{T}}\left( \nabla w\right) ^{2}\psi _{\lambda
}dxdt+ \\ 
+C_{1}\left( 1/\lambda ^{3}\right) \dint\limits_{Q_{T}}\left[ \left( \Delta
w\right) ^{2}+\left( \nabla w\right) ^{2}\right] \psi _{\lambda }dxdt+ \\ 
+2\dint\limits_{Q_{T}}\left\vert \chi _{\varepsilon }^{\prime }\left(
t\right) \right\vert ^{2}\left( \widetilde{v}^{2}+\widetilde{m}^{2}\right)
\psi _{\lambda }dxdt+C_{1}\delta ^{2}\lambda ^{2}e^{2\lambda b^{2}},\text{ }%
\forall \lambda \geq \lambda _{1}.%
\end{array}%
\right.  \label{4.28}
\end{equation}%
Apply Carleman estimate (\ref{3.4}) of Theorem 4.1 for the operator $%
\partial _{t}+\alpha \Delta $ to the left hand side of (\ref{4.28}) and note
that $1/\lambda ^{3}<<1/\lambda $ at $\lambda \rightarrow \infty .$ Also
note that (\ref{4.220}) and (\ref{4.221}) imply that the last line of (\ref%
{3.4}) equals zero if $u$ is replaced with $w.$ We obtain for sufficiently
large $\lambda _{1}:$%
\begin{equation}
\left. 
\begin{array}{c}
\left( 1/\lambda \right) \dint\limits_{Q_{T}}w_{t}^{2}\psi _{\lambda
}dxdt+\left( 1/\lambda \right)
\dsum\limits_{i,j=1}^{n}\dint\limits_{Q_{T}}w_{x_{i}x_{j}}^{2}\psi _{\lambda
}dxdt+ \\ 
+\lambda \dint\limits_{Q_{T}}\left( \nabla w\right) ^{2}\psi _{\lambda
}dxdt+\lambda ^{3}\dint\limits_{Q_{T}}w^{2}\psi _{\lambda }dxdt\leq \\ 
\leq 2\dint\limits_{Q_{T}}\left\vert \chi _{\varepsilon }^{\prime }\left(
t\right) \right\vert ^{2}\left( \widetilde{v}^{2}+\widetilde{m}^{2}\right)
\psi _{\lambda }dxdt+C_{1}\delta ^{2}\lambda ^{2}e^{2\lambda b^{2}},\text{ }%
\forall \lambda \geq \lambda _{1}.%
\end{array}%
\right.  \label{4.29}
\end{equation}%
By (\ref{4.29}) the\ term 
\begin{equation*}
C_{1}\dint\limits_{Q_{T}}\left[ \left( \Delta w\right) ^{2}+\left( \nabla
w\right) ^{2}\right] \psi _{\lambda }dxdt
\end{equation*}%
in the third line of (\ref{4.26}) can be estimated as:%
\begin{equation}
\left. 
\begin{array}{c}
C_{1}\dint\limits_{Q_{T}}\left[ \left( \Delta w\right) ^{2}+\left( \nabla
w\right) ^{2}\right] \psi _{\lambda }dxdt\leq \\ 
\leq 2\lambda \dint\limits_{Q_{T}}\left\vert \chi _{\varepsilon }^{\prime
}\left( t\right) \right\vert ^{2}\left( \widetilde{v}^{2}+\widetilde{m}%
^{2}\right) \psi _{\lambda }dxdt+C_{1}\delta ^{2}\lambda ^{2}e^{2\lambda
b^{2}}.%
\end{array}%
\right.  \label{4.290}
\end{equation}%
Hence, it follows from (\ref{4.26}), (\ref{4.29}) and (\ref{4.290}) that for
sufficiently large $\lambda _{1}$%
\begin{equation}
\left. 
\begin{array}{c}
\dint\limits_{Q_{T}}\left( q_{t}^{2}+w_{t}^{2}\right) \psi _{\lambda
}dxdt+\dsum\limits_{i,j=1}^{n}\dint\limits_{Q_{T}}\left(
q_{x_{i}x_{j}}^{2}+w_{x_{i}x_{j}}^{2}\right) \psi _{\lambda }dxdt+ \\ 
+\dint\limits_{Q_{T}}\left[ \left( \nabla q\right) ^{2}+\left( \nabla
w\right) ^{2}+\left( q^{2}+w^{2}\right) \right] \psi _{\lambda }dxdt\leq \\ 
\leq 2\lambda \dint\limits_{Q_{T}}\left\vert \chi _{\varepsilon }^{\prime
}\left( t\right) \right\vert ^{2}\left( \widetilde{v}^{2}+\widetilde{m}%
^{2}\right) \psi _{\lambda }dxdt+C_{1}\delta ^{2}\lambda ^{2}e^{2\lambda
b^{2}},\text{ }\forall \lambda \geq \lambda _{1}.%
\end{array}%
\right.  \label{4.30}
\end{equation}%
\newline

Consider the first term in the last line of (\ref{4.30}). By (\ref{2.5}) and
(\ref{4.220})%
\begin{equation}
\left. 
\begin{array}{c}
2\lambda \dint\limits_{Q_{T}}\left\vert \chi _{\varepsilon }^{\prime }\left(
t\right) \right\vert ^{2}\left( \widetilde{v}^{2}+\widetilde{m}^{2}\right)
\psi _{\lambda }dxdt= \\ 
=2\lambda \dint\limits_{0}^{\varepsilon }\left\vert \chi _{\varepsilon
}^{\prime }\left( t\right) \right\vert ^{2}\left[ \dint\limits_{\Omega
}\left( \widetilde{v}^{2}+\widetilde{m}^{2}\right) \psi _{\lambda }dx\right]
dt+ \\ 
+2\lambda \dint\limits_{T-\varepsilon }^{T}\left\vert \chi _{\varepsilon
}^{\prime }\left( t\right) \right\vert ^{2}\left[ \dint\limits_{\Omega
}\left( \widetilde{v}^{2}+\widetilde{m}^{2}\right) \psi _{\lambda }dx\right]
dt\leq \\ 
\leq C_{1}\lambda \left[ \dint\limits_{0}^{\varepsilon }\left(
\dint\limits_{\Omega }\psi _{\lambda }dx\right)
dt+\dint\limits_{T-\varepsilon }^{T}\left( \dint\limits_{\Omega }\psi
_{\lambda }dx\right) dt\right] \times \\ 
\times \left( \left\Vert \widetilde{v}\right\Vert _{L_{2}\left( Q_{T}\right)
}^{2}+\left\Vert \widetilde{m}\right\Vert _{L_{2}\left( Q_{T}\right)
}^{2}\right) ,\text{ }\forall \lambda >0.%
\end{array}%
\right.  \label{4.31}
\end{equation}%
It follows from (\ref{2.1}) and (\ref{3.1}) that 
\begin{equation*}
\psi _{\lambda }\left( x_{1},t\right) \leq \exp \left[ 2\lambda \left(
b^{2}-c^{2}\left( \frac{T}{2}-\varepsilon \right) ^{2}\right) \right] \text{
in }\Omega \times \left( \left( 0,\varepsilon \right) \cup \left(
T-\varepsilon ,T\right) \right) ,\text{ }\forall \lambda >0.
\end{equation*}%
Hence, (\ref{4.31}) implies%
\begin{equation}
\left. 
\begin{array}{c}
2\lambda \dint\limits_{0}^{\varepsilon }\left\vert \chi _{\varepsilon
}^{\prime }\left( t\right) \right\vert ^{2}\left[ \dint\limits_{\Omega
}\left( \widetilde{v}^{2}+\widetilde{m}^{2}\right) \psi _{\lambda }dx\right]
dt\leq \\ 
\leq C_{1}\lambda \exp \left[ 2\lambda \left( b^{2}-c^{2}\left(
T/2-\varepsilon \right) ^{2}\right) \right] \left( \left\Vert \widetilde{v}%
\right\Vert _{L_{2}\left( Q_{T}\right) }^{2}+\left\Vert \widetilde{m}%
\right\Vert _{L_{2}\left( Q_{T}\right) }^{2}\right) ,\text{ }\forall \lambda
>0.%
\end{array}%
\right.  \label{4.32}
\end{equation}

Consider now the domain $Q_{T,2\varepsilon },$ which was defined in (\ref%
{3.01}). By (\ref{3.1})%
\begin{equation}
\psi _{\lambda }\left( x_{1},t\right) \geq \exp \left[ 2\lambda \left(
a^{2}-c^{2}\left( T/2-2\varepsilon \right) ^{2}\right) \right] \text{ in }%
Q_{T,2\varepsilon }.  \label{4.33}
\end{equation}%
Since $Q_{T,2\varepsilon }\subset Q_{T},$ then (\ref{3.100}), (\ref{4.30}), (%
\ref{4.32}) and (\ref{4.33}) lead to the following estimate:%
\begin{equation}
\left. 
\begin{array}{c}
\exp \left[ 2\lambda \left( a^{2}-c^{2}\left( T/2-2\varepsilon \right)
^{2}\right) \right] \left[ \left\Vert \widetilde{v}\right\Vert
_{H^{2,1}\left( Q_{T,2\varepsilon }\right) }^{2}+\left\Vert \widetilde{m}%
\right\Vert _{H^{2,1}\left( Q_{T,2\varepsilon }\right) }^{2}\right] \leq \\ 
\leq C_{1}\delta ^{2}\lambda ^{2}e^{2\lambda b^{2}}+ \\ 
+C_{1}\lambda ^{2}\exp \left[ 2\lambda \left( b^{2}-c^{2}\left(
T/2-\varepsilon \right) ^{2}\right) \right] \left( \left\Vert \widetilde{v}%
\right\Vert _{L_{2}\left( Q_{T}\right) }^{2}+\left\Vert \widetilde{m}%
\right\Vert _{L_{2}\left( Q_{T}\right) }^{2}\right) , \\ 
\text{ }\forall \lambda \geq \lambda _{1}.%
\end{array}%
\right.  \label{4.330}
\end{equation}%
Dividing (\ref{4.330}) by $\exp \left[ 2\lambda \left( a^{2}-c^{2}\left(
T/2-2\varepsilon \right) ^{2}\right) \right] $, keeping in mind that by (\ref%
{4.02}) and (\ref{4.03}) $\varepsilon \in \left( 0,T/4\right) $ and using (%
\ref{4.142}), we obtain for all $\lambda \geq \lambda _{1}:$ 
\begin{equation}
\left. 
\begin{array}{c}
\left\Vert \widetilde{v}\right\Vert _{H^{2,1}\left( Q_{T,2\varepsilon
}\right) }^{2}+\left\Vert \widetilde{m}\right\Vert _{H^{2,1}\left(
Q_{T,2\varepsilon }\right) }^{2}\leq \\ 
\leq C_{1}\delta ^{2}\lambda ^{2}e^{2\lambda d}+ \\ 
+C_{1}\lambda ^{2}\exp \left[ -2\lambda \left( c^{2}\varepsilon \left(
T-3\varepsilon \right) -\left( b^{2}-a^{2}\right) \right) \right] \left(
\left\Vert \widetilde{v}\right\Vert _{L_{2}\left( Q_{T}\right)
}^{2}+\left\Vert \widetilde{m}\right\Vert _{L_{2}\left( Q_{T}\right)
}^{2}\right) ,\text{ }%
\end{array}%
\right.  \label{4.34}
\end{equation}%
\begin{equation}
d=\left( b^{2}-a^{2}\right) \left[ 1+\frac{\left( 1+\beta \right) \left(
T/2-2\varepsilon \right) ^{2}}{\varepsilon \left( T-3\varepsilon \right) }%
\right] .  \label{4.340}
\end{equation}%
Using (\ref{4.142}), we obtain%
\begin{equation*}
-2\lambda \left( c^{2}\varepsilon \left( T-3\varepsilon \right) -\left(
b^{2}-a^{2}\right) \right) =-2\lambda \beta \left( b^{2}-a^{2}\right) .
\end{equation*}%
Hence, (\ref{4.34}) and (\ref{4.340}) imply 
\begin{equation}
\left. 
\begin{array}{c}
\left\Vert \widetilde{v}\right\Vert _{H^{2,1}\left( Q_{T,2\varepsilon
}\right) }+\left\Vert \widetilde{m}\right\Vert _{H^{2,1}\left(
Q_{T,2\varepsilon }\right) }\leq C_{1}\delta \lambda e^{\lambda d}+ \\ 
+C_{1}\lambda e^{-\lambda \beta \left( b^{2}-a^{2}\right) }\left( \left\Vert 
\widetilde{v}\right\Vert _{L_{2}\left( Q_{T}\right) }+\left\Vert \widetilde{m%
}\right\Vert _{L_{2}\left( Q_{T}\right) }\right) ,\text{ }\forall \lambda
\geq \lambda _{1}.%
\end{array}%
\right.  \label{4.35}
\end{equation}%
We choose $\lambda =\lambda \left( \delta \right) $ such that 
\begin{equation}
e^{-\lambda \beta \left( b^{2}-a^{2}\right) }=\delta e^{\lambda d}.
\label{4.350}
\end{equation}%
Hence, 
\begin{equation}
\lambda =\lambda \left( \delta \right) =\ln \left( \delta ^{-1/\left(
d+\beta \left( b^{2}-a^{2}\right) \right) }\right) .  \label{4.36}
\end{equation}%
The choice (\ref{4.36}) is possible if $\delta \in \left( 0,\delta
_{0}\right) $ if the number $\delta _{0}=\delta _{0}\left( \alpha ,\rho
,\varepsilon ,R,Q_{T}\right) \in \left( 0,1\right) $ is so small that $\ln
\left( \delta _{0}^{-1/\left( d+\beta \left( b^{2}-a^{2}\right) \right)
}\right) \geq \lambda _{1}.$ Hence, (\ref{4.35}) becomes:%
\begin{equation}
\left. 
\begin{array}{c}
\left\Vert \widetilde{v}\right\Vert _{H^{2,1}\left( Q_{T,2\varepsilon
}\right) }+\left\Vert \widetilde{m}\right\Vert _{H^{2,1}\left(
Q_{T,2\varepsilon }\right) }\leq C_{1}\lambda \left( \delta \right) \delta
^{\rho }\left( 1+\left\Vert \widetilde{v}\right\Vert _{L_{2}\left(
Q_{T}\right) }+\left\Vert \widetilde{m}\right\Vert _{L_{2}\left(
Q_{T}\right) }\right) , \\ 
\text{ }\forall \delta \in \left( 0,\delta _{0}\right) ,%
\end{array}%
\right.  \label{4.37}
\end{equation}%
\begin{equation}
\rho =\frac{\beta \left( b^{2}-a^{2}\right) }{d+\beta \left(
b^{2}-a^{2}\right) }.  \label{4.370}
\end{equation}

We now show that for any value of the number $\rho $ as in (\ref{4.01})
there exists unique number $\beta =\beta \left( \rho ,\varepsilon ,T\right)
>0$ such that (\ref{4.370}) holds. Recall that the positivity of the
parameter $\beta $ is required in (\ref{4.142}). By (\ref{4.370}) 
\begin{equation}
\beta \left( b^{2}-a^{2}\right) \left( 1-\rho \right) =\rho d.  \label{4.38}
\end{equation}%
Substituting in (\ref{4.38}) the formula for $d$ from (\ref{4.340}), we
obtain the following equation for $\beta :$%
\begin{equation}
\beta \frac{\left( 1-\rho \right) \varepsilon \left( T-3\varepsilon \right)
-\rho \left( T/2-2\varepsilon \right) ^{2}}{\varepsilon \left(
T-3\varepsilon \right) }=\rho \left[ 1+\frac{\left( T/2-2\varepsilon \right)
^{2}}{\varepsilon \left( T-3\varepsilon \right) }\right] .  \label{4.39}
\end{equation}%
Recall that by (\ref{4.03}) $T-3\varepsilon >0.$ Hence, the solution $\beta $
of equation (\ref{4.39}) is positive if and only if 
\begin{equation*}
\left( 1-\rho \right) \varepsilon \left( T-3\varepsilon \right) -\rho \left(
T/2-2\varepsilon \right) ^{2}>0,
\end{equation*}%
which is equivalent with%
\begin{equation*}
\left( \rho +3\right) \varepsilon ^{2}-\varepsilon T+\rho \frac{T^{2}}{4}<0.
\end{equation*}%
Hence, we must have 
\begin{equation*}
\frac{T}{2\left( \rho +3\right) }\left( 1-\sqrt{1-3\rho -\rho ^{2}}\right)
<\varepsilon <\frac{T}{2\left( \rho +3\right) }\left( 1+\sqrt{1-3\rho -\rho
^{2}}\right) .
\end{equation*}%
Hence, (\ref{4.01})-(\ref{4.03}) imply that the solution $\beta $ of
equation (\ref{4.39}) is indeed positive for the values of $\rho $
satisfying (\ref{4.01}).

Next, having this $\beta ,$ consider an arbitrary number $\gamma \in \left(
0,\rho \right) .$ It follows from (\ref{4.36}) that there exists a
sufficiently small number $\delta _{1}=\delta _{1}\left( \alpha ,\rho
,\gamma ,\varepsilon ,R,Q_{T}\right) \in \left( 0,\delta _{0}\left( \rho
,\varepsilon ,R,Q_{T}\right) \right) $ such that 
\begin{equation*}
\lambda \left( \delta \right) =\ln \left( \delta ^{-1/\left( d+\beta \left(
b^{2}-a^{2}\right) \right) }\right) <\delta ^{-\gamma },\text{ }\forall
\delta \in \left( 0,\delta _{1}\right) .
\end{equation*}%
Substituting this in (\ref{4.37}) and keeping in mind notations in the first
line of (\ref{4.16}), we arrive at the target estimate (\ref{4.13}) of this
theorem.

We now prove the uniqueness statement of this theorem. Suppose that
conditions (\ref{4.140}) hold. This means that we can set $\delta =0$ in (%
\ref{4.13}), which implies that $\widetilde{m}\left( x,t\right) =\widetilde{v%
}\left( x,t\right) =0$ in $Q_{T,2\varepsilon }.$ Setting $\rho \rightarrow 0$
in the last line of (\ref{4.03}), we obtain that we can take an arbitrary
small $\varepsilon >0.$ Hence, $\widetilde{m}\left( x,t\right) \equiv 
\widetilde{v}\left( x,t\right) \equiv 0$ in $Q_{T}.$ $\square $

\subsection{The case of a sufficiently large observation time $T>0$}

\label{sec:5.2}

To avoid working with cumbersome formulas and, thus, to simplify the
presentation this way, we specify here the dependence of $\varepsilon $ from 
$T$ in (\ref{3.01}). Thus, (\ref{3.01}) is modified now as:%
\begin{equation}
\left. 
\begin{array}{c}
\varepsilon =T/6<T/4, \\ 
Q_{T,T/6}=\Omega \times \left( T/6,5T/6\right) , \\ 
Q_{T,T/3}=\Omega \times \left( T/3,2T/3\right) \subset Q_{T,T/6}\subset
Q_{T}.%
\end{array}%
\right.  \label{4.40}
\end{equation}%
Also, we assign the parameter $c\equiv 1$ in CWF (\ref{3.1}). Hence, now (%
\ref{3.1}) and (\ref{3.100})\ became%
\begin{equation}
\left. 
\begin{array}{c}
\psi _{\lambda }\left( x_{1},t\right) =\exp \left[ \lambda \left(
x_{1}^{2}-\left( t-T/2\right) ^{2}\right) \right] , \\ 
\exp \left[ 2\lambda \left( a^{2}-\left( t-T/2\right) ^{2}\right) \right]
\leq \psi _{\lambda }^{2}\left( x_{1},t\right) \leq e^{2\lambda b^{2}}\text{
in }Q_{T}.%
\end{array}%
\right.  \label{4.41}
\end{equation}

\textbf{Theorem 5.2} (the case of a sufficiently large $T).$ \emph{Suppose
that conditions (\ref{4.40}) hold. Also, let assumptions of items number 2-5
of \ Theorem 5.1 be valid. Suppose that }%
\begin{equation}
\left( b^{2}-a^{2}\right) =\frac{T^{2}}{12}\left( 1-s\right) ,\text{ }s\in
\left( 0,1\right) .  \label{4.42}
\end{equation}%
\emph{Choose an arbitrary number }$\sigma $\emph{\ such that} 
\begin{equation}
\sigma \in \left( 0,\frac{3}{4}s\right) .  \label{4.420}
\end{equation}%
\emph{Then there exists a sufficiently small number }$\delta _{2}=\delta
_{2}\left( \alpha ,\sigma ,s,R,Q_{T}\right) \in \left( 0,1\right) $\emph{\
such that the following analog of the H\"{o}lder stability estimate (\ref%
{4.13}) of Theorem 5.1 holds:}%
\begin{equation}
\left. 
\begin{array}{c}
\left\Vert u_{1}-u_{2}\right\Vert _{H^{2,1}\left( Q_{T,T/3}\right)
}+\left\Vert m_{1}-m_{2}\right\Vert _{H^{2,1}\left( Q_{T,T/3}\right) }\leq
\\ 
\leq C_{2}\left( 1+\left\Vert u_{1}-u_{2}\right\Vert _{L_{2}\left(
Q_{T}\right) }+\left\Vert m_{1}-m_{2}\right\Vert _{H^{2}\left( Q_{T}\right)
}\right) \delta ^{3s/4-\sigma },\text{ }\forall \delta \in \left( 0,\delta
_{2}\right) ,\text{ }%
\end{array}%
\right.  \label{4.43}
\end{equation}%
\emph{where the numbers }$\delta _{2}\left( \alpha ,\sigma ,s,R,Q_{T}\right) 
$ \emph{and }$C_{2}=C_{2}\left( \alpha ,R,Q_{T}\right) >0$\emph{\ depend
only on listed parameters.}

\textbf{Remarks 5.2.}

\begin{enumerate}
\item \emph{Estimate (\ref{4.43}) is stronger than estimate (\ref{4.13})
since (\ref{4.01}) implies that }$\rho <0.31<0.75s=3s/4$ for $s\in \left(
0.42,1\right) $ \emph{and }$\sigma \in \left( 0,3s/4\right) $\emph{\ is an
arbitrary number in (\ref{4.420}). This is quite natural from the intuition
standpoint since the larger the observation time is, the more informative
the data are, and therefore, the more stable the problem should be. In our
case, these are the lateral Cauchy data (\ref{4.9}), (\ref{4.10}).}

\item \emph{It follows immediately from Theorem 5.1 that uniqueness of the
solution holds in the case of conditions of Theorem 5.2 as well. Hence, we
do not include the claim of uniqueness in the formulation of Theorem 5.2.}
\end{enumerate}

\textbf{Proof of Theorem 5.2.} In this proof $C_{2}=C_{2}\left( \alpha
,R,Q_{T}\right) >0$\emph{\ \ }denotes different constants depending only on
listed parameters.

We keep the same notations as the ones in the proof of Theorem 5.1. Until (%
\ref{4.330}) the steps of this proof are the same as the ones in the proof
of Theorem 5.1. \ Hence, we do not repeat those steps now. Using (\ref{4.40}%
) and (\ref{4.41}), we obtain the following direct analog of (\ref{4.330}):

\begin{equation}
\left. 
\begin{array}{c}
\exp \left[ 2\lambda \left( a^{2}-T^{2}/36\right) \right] \left[ \left\Vert 
\widetilde{v}\right\Vert _{H^{2,1}\left( Q_{T,T/3}\right) }^{2}+\left\Vert 
\widetilde{m}\right\Vert _{H^{2,1}\left( Q_{T,T/3}\right) }^{2}\right] \leq
\\ 
\leq C_{1}\delta ^{2}\lambda ^{2}e^{2\lambda b^{2}}+ \\ 
+C_{1}\lambda ^{2}\exp \left[ 2\lambda \left( b^{2}-T^{2}/9\right) \right]
\left( \left\Vert \widetilde{v}\right\Vert _{L_{2}\left( Q_{T}\right)
}^{2}+\left\Vert \widetilde{m}\right\Vert _{L_{2}\left( Q_{T}\right)
}^{2}\right) ,\text{ }\forall \lambda \geq \lambda _{2},\text{ }%
\end{array}%
\right.  \label{4.50}
\end{equation}%
where $\lambda _{2}=\lambda _{2}\left( \alpha ,R,Q_{T}\right) \geq \lambda
_{0}$ is the analog of the parameter $\lambda _{1}.$ Dividing both parts of (%
\ref{4.50}) by $\exp \left[ 2\lambda \left( a^{2}-T^{2}/36\right) \right] ,$
we obtain%
\begin{equation}
\left. 
\begin{array}{c}
\left\Vert \widetilde{v}\right\Vert _{H^{2,1}\left( Q_{T,T/3}\right)
}+\left\Vert \widetilde{m}\right\Vert _{H^{2,1}\left( Q_{T,T/3}\right) }\leq
\\ 
\leq C_{2}\delta \lambda \exp \left[ \lambda \left(
b^{2}-a^{2}+T^{2}/36\right) \right] + \\ 
+C_{2}\lambda \exp \left[ -\lambda \left( T^{2}/12-\left( b^{2}-a^{2}\right)
\right) \right] \left( \left\Vert \widetilde{v}\right\Vert _{L_{2}\left(
Q_{T}\right) }+\left\Vert \widetilde{m}\right\Vert _{L_{2}\left(
Q_{T}\right) }\right) ,\text{ }\forall \lambda \geq \lambda _{2}.%
\end{array}%
\right.  \label{4.52}
\end{equation}

Similarly with (\ref{4.350}) we choose $\lambda =\lambda \left( \delta
\right) $ such that 
\begin{equation}
\exp \left[ -\lambda \left( T^{2}/12-\left( b^{2}-a^{2}\right) \right) %
\right] =\delta \exp \left[ \lambda \left( b^{2}-a^{2}+T^{2}/36\right) %
\right] .  \label{4.53}
\end{equation}%
Hence,%
\begin{equation}
\lambda =\lambda \left( \delta \right) =\ln \left( \delta ^{-9/T^{2}}\right)
.  \label{4.54}
\end{equation}%
Similarly with (\ref{4.36}) the choice (\ref{4.54}) is possible if $\delta
\in \left( 0,\widehat{\delta }\right) $, where the number $\widehat{\delta }=%
\widehat{\delta }\left( \alpha ,R,Q_{T}\right) \in \left( 0,1\right) $ is so
small that $\ln \left( \widehat{\delta }^{-9/T^{2}}\right) \geq \lambda _{2}.
$ By (\ref{4.42})%
\begin{equation*}
\exp \left[ -\lambda \left( T^{2}/12-\left( b^{2}-a^{2}\right) \right) %
\right] =\exp \left( -\lambda sT^{2}/12\right) .
\end{equation*}%
Hence, by (\ref{4.54})  
\begin{equation}
\exp \left[ -\lambda \left( \delta \right) \left( T^{2}/12-\left(
b^{2}-a^{2}\right) \right) \right] =\exp \left( -\lambda \left( \delta
\right) sT^{2}/12\right) =\delta ^{3s/4},\forall \delta \in \left( 0,%
\widehat{\delta }\right) .  \label{4.55}
\end{equation}%
Hence, (\ref{4.52})-(\ref{4.550}) imply%
\begin{equation}
\left. 
\begin{array}{c}
\left\Vert \widetilde{v}\right\Vert _{H^{2,1}\left( Q_{T,T/3}\right)
}+\left\Vert \widetilde{m}\right\Vert _{H^{2,1}\left( Q_{T,T/3}\right) }\leq 
\\ 
\leq C_{2}\delta ^{3s/4}\ln \left( \ln \left( \delta ^{-9/T^{2}}\right)
\right) \left( 1+\left\Vert \widetilde{v}\right\Vert _{L_{2}\left(
Q_{T}\right) }+\left\Vert \widetilde{m}\right\Vert _{L_{2}\left(
Q_{T}\right) }\right) ,\text{ }\forall \delta \in \left( 0,\widehat{\delta }%
\right) .%
\end{array}%
\right.   \label{1}
\end{equation}%
Next, by (\ref{4.420}) and (\ref{4.54}) there exists a sufficiently small
number $\delta _{2}=\delta _{2}\left( \alpha ,\sigma ,s,R,Q_{T}\right) \in
\left( 0,\widehat{\delta }\right) $ such that 
\begin{equation}
\ln \lambda \left( \delta \right) =\ln \left( \ln \left( \delta
^{-9/T^{2}}\right) \right) <\delta ^{-\sigma },\text{ }\forall \delta \in
\left( 0,\delta _{2}\right) .  \label{4.550}
\end{equation}%
It follows from (\ref{1}) and (\ref{4.550}) that%
\begin{equation}
\left. 
\begin{array}{c}
\left\Vert \widetilde{v}\right\Vert _{H^{2,1}\left( Q_{T,T/3}\right)
}+\left\Vert \widetilde{m}\right\Vert _{H^{2,1}\left( Q_{T,T/3}\right) }\leq 
\\ 
\leq C_{2}\delta ^{3s/4-\sigma }\left( 1+\left\Vert \widetilde{v}\right\Vert
_{L_{2}\left( Q_{T}\right) }+\left\Vert \widetilde{m}\right\Vert
_{L_{2}\left( Q_{T}\right) }\right) ,\text{ }\forall \delta \in \left(
0,\delta _{2}\right) .%
\end{array}%
\right.   \label{4.56}
\end{equation}%
The target estimate (\ref{4.43}) of this theorem follows immediately from
notations in the first line of (\ref{4.16}) and estimate (\ref{4.56}). $\
\square $


\begin{thebibliography}{99}
\bibitem{A} Y. Achdou, P. Cardaliaguet, F. Delarue, A. Porretta and F.
Santambrogio, \emph{Mean Field Games}, Cetraro, Italy 2019, Lecture Notes in
Mathematics, C.I.M.E. Foundation Subseries, Volume 2281, Springer, 2019.

\bibitem{A1} Y. Achdou and J.-M. Lasry, Mean field games for modeling crowd
motion. \emph{Contributions to Partial Differential Equations and
Applications, }17--42,\emph{\ Comput. Methods Appl. Sci., }47, Springer,
Cham, 2019.\emph{\ }

\bibitem{A2} D.M. Ambrose, Existence theory for a time--dependent mean field
games model of household wealth,~\emph{Appl. Math. Optim.,}~83, 2051--2081,
2021.

\bibitem{B} D. Bauso, H. Tembine and T. Basar, Opinion dynamics in social
networks through mean-field games, \emph{SIAM J. Control Optim.,} 54,
3225--3257, 2016.

\bibitem{BukhKlib} A.~Bukhgeim and M.~V. Klibanov, Uniqueness in the large
of a class of multidimensional inverse problems, \emph{Soviet Mathematics
Doklady}, 17, 244--247, 1981.

\bibitem{Burger} M.~Burger, L.~Caffarelli, and P.~A. Markowich, Partial
differential equation models in the socio-economic sciences, Philosophical
Transactions of Royal Society, A372, 20130406, 2014.

\bibitem{Chow} Y.T. Chow, S.W. Fung, S. Liu, L. Nurbekyan, and S. Osher, A
numerical algorithm for inverse problem from partial boundary measurement
arising from mean field game problem, \emph{Inverse Problems}, 39, 014001,
2023.

\bibitem{Cou} R. Couillet, S.M. Perlaza, H. Tembine, and M. Debbah,
Electrical vehicles in the smart grid: a mean field game analysis, \emph{%
IEEE Journal on Selected Areas of Communications}, 30, 1086-1096, 2012.

\bibitem{Huang1} M.~Huang, P.~E. Caines, and R.~P. Malham\'{e},
Large-population cost-coupled LQG problems with nonuniform agents:
individual-mass behavior and decentralized Nash equilibria, \emph{IEEE
Trans. Automat. Control}, 52,~1560--1571, 2007.

\bibitem{Huang2} M.~Huang, R.~P. Malham\'{e} and P.~E. Caines, Large
population stochastic dynamic games: closed-loop McKean-Vlasov systems and
the Nash certainty equivalence principle, \emph{Commun. Inf. Syst.}, 6,
221--251, 2006.

\bibitem{Ksurvey} M.V. Klibanov, Carleman estimates for global uniqueness,
stability and numerical methods for coefficient inverse problems, \emph{J.
of Inverse and Ill-Posed Problems}, 21, 477-510, 2013.

\bibitem{KLpar} {M.V. Klibanov, J. Li and W. Zhang}, Convexification for an
inverse parabolic problem, \emph{Inverse Problems}, 36: 085008, 2020.

\bibitem{KL} M.V. Klibanov and J. Li, \emph{Inverse Problems and Carleman
Estimates: Global Uniqueness, Global Convergence and Experimental Data}, De
Gruyter, 2021.

\bibitem{MFG1} M.V. Klibanov and Y. Averboukh, Lipschitz stability estimate
and uniqueness in the retrospective analysis for the mean field games system
via two Carleman estimates, \emph{arXiv}: 2302.10709v2, accepted in \emph{%
SIAM J. Mathematical Analysis, }2023.

\bibitem{MFG2} M.V. Klibanov, The mean field games system: Carleman
estimates, Lipschitz stability and uniqueness, \emph{J. of Inverse and
Ill-Posed Problems}, published online,
https://doi.org/10.1515/jiip-2023-0023, 2023.

\bibitem{MFG4} M.V. Klibanov, J. Li and H. Liu, H\"{o}lder stability and
uniqueness for the mean field games system via Carleman estimates, \emph{%
Studies in Applied Mathematics}, https://doi.org/10.1111/sapm.12633, 1-24,
2023.

\bibitem{MFG5} M.~V. Klibanov, J.~Li, and H.~Liu, Coefficient inverse
problems for a generalized mean field games system with the final
overdetermination, \emph{arXiv}: 2305.01065, 2023.

\bibitem{MFG6} M.V. Klibanov, A coefficient inverse problem for the mean
field games system, \emph{Applied Mathematics and Optimization,} 88:54, 2023.

\bibitem{MFG7} M.V. Klibanov, J. Li and Z. Yang, Convexification numerical
method for the retrospective problem of mean field games, \emph{arXiv}:
2306.14404, 2023.

\bibitem{MFG8} M.~V. Klibanov, J.~Li, and Z.~Yang. Convexification for a
coefficient inverse problem of mean field games, \emph{arXiv}: 2310.08878,
2023.

\bibitem{KB} V.N. Kolokoltsov and A. Bensoussan, Mean-field-game model of
botnet defence in cybersecurity, \emph{Applied Mathematics and Optimization}%
, 74, 669--692, 2015.

\bibitem{KM} V.N. Kolokoltsov, O.A. Malafeyev, Mean field game model of
corruption, \emph{Dynamics Games and Applications} 7, 34--47, 2017.

\bibitem{Kol} V.N. Kolokoltsov and O. A. Malafeyev, \emph{Many Agent Games
in Socio-economic Systems: Corruption, Inspection, Coalition Building,
Network Growth, Security}, Springer Nature Switzerland AG, 2019.

\bibitem{LL1} J.-M. Lasry and P.-L. Lions, Jeux \`{a} champ moyen. I. Le cas
stationnaire, \emph{C. R. Math. Acad. Sci. Paris,} 343, 619--625, 2006.

\bibitem{LL2} J.-M. Lasry and P.-L. Lions, Jeux \`{a} champ moyen. II.
Horizon fini et contr\^{o}le optimal, \emph{C. R. Math. Acad. Sci. Paris},
343, 679--684, 2006.

\bibitem{LL3} J.-M. Lasry and P.-L. Lions, Mean field games, \emph{Japanese
Journal of Mathematics}, 2, 229-260, 2007.

\bibitem{LRS} M.M. Lavrentiev, V.G. Romanov and S.P. Shishatskii, \emph{%
Ill-Posed Problems of Mathematical Physics and Analysis}, AMS, Providence,
RI, 1986.

\bibitem{Liu1} H.~Liu and S.~Zhang, On an inverse boundary problem for mean
field games, \emph{arXiv}: 2212.09110, 2022.

\bibitem{Liu2} H.~Liu and S.~Zhang. Simultaneously recovering running cost
and Hamiltonian in mean field games system, \emph{arXiv}:2303.13096, 2023.

\bibitem{Liu3} H.~Liu, C.~Mou, and S.~Zhang, Inverse problems for mean field
games, \emph{Inverse Problems}, 39, 085003, 2023.

\bibitem{Nov1} R. G. Novikov, The inverse scattering problem on a fixed
energy level for the two-dimensional Schr\"{o}dinger operator, \emph{J.
Functional Analysis,} 103, 409-463, 1992.

\bibitem{Nov2} R. G. Novikov, $\partial -$bar approach to approximate
inverse scattering at fixed energy in three dimensions, \emph{International
Math. Research Peports,} 6,\textbf{\ }287-349, 2005.

\bibitem{Rom} V. G. Romanov, \emph{Inverse Problems of Mathematical Physics}%
, VNU Press, Utrecht, 1987.

\bibitem{T} A. N. Tikhonov and V.Ya. Arsenin, \emph{Solutions of Ill-Posed
Problems}, Winston\&Sons, Washington, \ D.C., 1977.

\bibitem{Trusov} N.~V. Trusov, Numerical study of the stock market crises
based on mean field games approach, \emph{Journal of Inverse and Ill-posed
Problems,} 29 (2021), pp.~849--865.

\bibitem{V} M. M. Vajnberg, \emph{Variational Method and Method of Monotone
Operators in the Theory of Nonlinear Equations}, Israel Program for
Scientific Translations, Jerusalem, 1973.

\bibitem{Yam} M. Yamamoto, Carleman estimates for parabolic equations\emph{.}
Topical Review, \emph{Inverse Problems,} 25\textbf{: }123013, 2009.
\end{thebibliography}
\end{document}